\def\isOneColumn{false}

\ifnum\pdfstrcmp{\isOneColumn}{false}=0
\documentclass{IEEEtran}
\usepackage{setspace}
\fi

\usepackage[noadjust]{cite}

\usepackage[caption=false,font=footnotesize]{subfig}

\usepackage{mathptmx} 
\usepackage{amsmath} 
\usepackage{amssymb}  
\usepackage{color}

\usepackage{accents}
\newcommand{\ubar}[1]{\underaccent{\bar}{#1}}

\newcommand{\norm}[1]{\lVert #1 \rVert}

\newcommand{\abs}[1]{\lvert #1 \rvert}

\DeclareMathOperator{\diag}{diag}
\DeclareMathOperator*{\minimize}{minimize}
\DeclareMathOperator*{\subjectto}{subject\ to}
\DeclareMathOperator{\tr}{tr}
\DeclareMathOperator*{\esssup}{ess\,sup}

\DeclareSymbolFont{bbold}{U}{bbold}{m}{n}
\DeclareSymbolFontAlphabet{\mathbbold}{bbold}

\usepackage{graphicx}

\newcommand{\onev}{\mathbbold{1}}

\usepackage{scalerel}
\newcommand{\nsubtheta}{{n_{\scaleto{\theta}{3pt}}}}

\usepackage{mathtools}

\usepackage{bm}

\newtheorem{definition}{Definition}[section]
\newtheorem{assumption}[definition]{Assumption}
\newtheorem{lemma}[definition]{Lemma}
\newtheorem{proposition}[definition]{Proposition}
\newtheorem{theorem}[definition]{Theorem}
\newtheorem{remark}[definition]{Remark}
\newtheorem{corollary}[definition]{Corollary}

\newtheorem{problem}[definition]{Problem}

\usepackage{calrsfs}
\DeclareMathAlphabet{\pazocal}{OMS}{zplm}{m}{n}
\renewcommand{\mathcal}[1]{\pazocal{#1}}

\newenvironment{proofof}[1]{
\renewcommand\IEEEproof{\noindent\hspace{2em}{\itshape Proof of #1: }}%
\begin{IEEEproof}}{\end{IEEEproof}
}

\usepackage{mathtools}
\mathtoolsset{centercolon}

\definecolor{forestgreen}{rgb}{0.33,0.61,0.34}

\newcommand{\afterequation}{\vskip 3pt}

\begin{document}

\title{Geometric Programming for Optimal Positive Linear Systems}
    
\author{Masaki Ogura, \IEEEmembership{Member, IEEE}, Masako Kishida, \IEEEmembership{Senior Member, IEEE}, and James Lam, \IEEEmembership{Fellow, IEEE}
\thanks{{The research was partly supported by ROIS NII Open Collaborative Research
2018 and GRF HKU {17201219}.}}
\thanks{M.~Ogura is with the Graduate School of Information Science and Technology, Osaka University, Suita, Osaka 565-0781, Japan. e-mail: {m-ogura@ist.osaka-u.ac.jp}}%
\thanks{M.~Kishida is with Principles of Informatics Research Division, National
Institute of Informatics, Tokyo 101-8430, Japan. e-mail: {kishida@nii.ac.jp}}%
\thanks{J.~Lam is with the Department of Mechanical Engineering, The
University of Hong Kong, Hong Kong. e-mail: {james.lam@hku.hk}}%
}

\maketitle

\begin{abstract}
This paper studies the parameter tuning problem of positive linear systems for optimizing their stability properties. We specifically show that, under certain regularity assumptions on the parametrization, the problem of finding the minimum-cost parameters that achieve a given requirement on a system norm reduces to a \emph{geometric program}, which in turn can be exactly and efficiently solved by convex optimization. The flexibility of geometric programming allows the state, input, and output matrices of the system to simultaneously depend on the parameters to be tuned. The class of system norms under consideration includes the $H^2$ norm, $H^\infty$ norm, Hankel norm, and Schatten $p$-norm. Also, the parameter tuning problem for ensuring the robust stability of the system under structural uncertainties is shown to be solved by geometric programming. The proposed optimization framework is further extended to delayed positive linear systems, where it is shown that the parameter tunning problem jointly constrained by the exponential decay rate, the $\mathcal L^1$-gain, and the $\mathcal L^\infty$\nobreakdash-gain can be solved by convex optimization. The assumption on the system parametrization is stated in terms of posynomial functions, which form a broad class of functions and thus allow us to deal with various interesting positive linear systems arising from, for example, dynamical buffer networks and epidemic spreading processes. We present numerical examples to illustrate the effectiveness of the proposed optimization framework. 
\end{abstract}


\section{Introduction}

Positive systems refer to, roughly speaking, the class of dynamical systems whose response signals to nonnegative input signals are constrained to be nonnegative~\cite{Farina2000,Rantzer2018}. The application areas in which positive systems naturally arise include  pharmacology~\cite{Hernandez-Vargas2013,Jonsson2014,Giordano2016a}, epidemiology~\cite{Nowzari2015a,Ogura2015c}, population biology~\cite{Belete2015,Ogura2016n}, 
and communication networks~\cite{Shorten2006}. In this context, several important advances towards the analysis and control of positive systems have been made in the last decade. For example, the authors in~\cite{Shorten2009} showed that stability of a positive linear system and the existence of a diagonal Lyapunov function are equivalent. It was shown in~\cite{HuijunGao2005} and~\cite{Rami2007} that structured stabilization problems for positive linear systems can be efficiently solved by linear matrix inequalities and a linear program, respectively. The authors in~\cite{Tanaka2011} showed that the celebrated Kalman-Yakubovich-Popov lemma admits a significantly simple representation in terms of diagonal quadratic storage functions.

Besides the aforementioned results concerning static-gain state-feedback control of positive linear systems, it has been observed in the literature that synthesis problems for positive linear systems often exhibit interesting convexity properties. For example, it was shown in~\cite{Colaneri2014a} that positive linear forms on the state variables of a time-varying positive linear system are convex with respect to the diagonals of its state matrix. The authors in~\cite{Dhingra2015} established the convexity of a symmetric modification of a class of steady-state disturbance attenuation problems. The authors in~\cite{Colombino2016b} showed the convexity of the power norm of output signals with respect to the diagonals of the state matrix. The authors in~\cite{Dhingra2017} presented an intrinsic convexity property of~$H^2$ and~$H^\infty$ state-feedback control problems for positive linear systems. A similar result is obtained in~\cite{Colombino2016} for robust state-feedback stabilization under structured uncertainties. However, the practical applicability of the aforementioned results is not necessarily enough to cover the wide range of applications of positive linear systems because the convexity properties in these results are mostly with respect to the diagonals of the state matrix of the system.

In this paper, we develop computationally efficient frameworks for tuning the parameters of a positive linear system, in which any entry of any of the state, input, and output matrices are allowed to be dependent on the parameter to be synthesized. We specifically show that, under certain regularity conditions on the parameterizations of these coefficient matrices, {the optimal parameter tuning problems constrained by the $H^2$ norm, $H^\infty$ norm, Hankel norm, and Schatten $p$-norm (for an even $p$) can be solved by \emph{geometric programming}~\cite{Boyd2007}. We also show that the problem of tuning the parameters for ensuring the robust stability of the system under structural uncertainties can be solved by geometric programming. We furthermore extend our framework to show that a class of mixed-constraint optimization problems for delayed positive linear systems can be solved by convex optimization.} A geometric program is a nonlinear optimization problem in which all the variables are positive and the objective function and constraints are described by monomial and posynomial functions (see Section~\ref{sec:prb} for details). Due to the log-log convexity of monomial and posynomial functions, a geometric program can be easily converted to an equivalent convex optimization problem, whose optimal solution can be efficiently found. {Furthermore, packages for directly formulating and solving geometric programs are available in various standard softwares including MATLAB, Python, and MOSEK. }As an illustration of our theoretical results, we study the buffer network optimization problem with $H^\infty$ norm constraints, and also the optimal medical resource allocation problem for robustly eradicating epidemic spreading processes taking place over uncertain complex networks~\cite{Preciado2014,Holme2017}.

Geometric programming has been successfully applied in various engineering areas including digital circuit design~\cite{Sapatnekar2004,Boyd2005}, chemical engineering~\cite{Wall1986}, power control in wireless networks~\cite{Chiang2007}, information theory~\cite{Chiang2004}, and structural design~\cite{Adeli1986} (see~\cite{Boyd2007} for an extensive list of applications). Since geometric programming offers a powerful tool for optimally tuning positive parameters, it would be natural to expect that this optimization framework allows us to synthesize positive systems as well. Despite this expectation, we find in the literature relatively few works for utilizing geometric programming to the synthesis of positive systems. An exception is the sequence of works~\cite{Preciado2014,Ogura2015a,Ogura2015i}, in which the authors study resource allocation problems for maximizing the exponential decay rate of the infection size within a networked epidemic spreading model. Although it was shown in~\cite{Ogura2015a} that a class of~$\mathcal L^1$-gain optimization problem for networked positive linear systems can be solved by geometric programming, it was not fully discussed in the reference if geometric programming applies to other classes of synthesis problems. It is finally remarked that other applications of geometric programming in the context of systems and control theory can be found in~\cite{Yazarel2004a,Singh2004}.

In this paper, we use the following notations. Let $\mathbb{R}$, $\mathbb{R}_+$, and~$\mathbb{R}_{++}$ denote the set of real, nonnegative, and positive numbers, respectively. For a positive integer~$n$, let $\{e_1, \dotsc, e_{n}\}$ denote the canonical basis of~$\mathbb{R}^{n}$. We let $\onev$ denote a column vector with all entries equal to one. The identity and the zero matrix of order~$n$ is denoted by~$I_n$ and $O_n$, respectively. A real matrix~$A$ is said to be nonnegative (positive), denoted by~$A\geq 0$ ($A>0$), if all entries of~$A$ are nonnegative (positive, respectively). We write $A\leq B$ if $B-A\geq 0$. The notations $A < B$, $A\geq B$, and~$A>B$ should be understood in the same manner. The maximum singular value of~$A$ is denoted by~$\norm{A}$. Let $A$ be a real and square matrix. We {say } that $A$ is Hurwitz stable if the eigenvalues of~$A$ have negative real parts. We say that $A$ is Metzler if the off-diagonal entries of~$A$ are nonnegative. By the Perron-Frobenius theorem~\cite{Horn1990}, a Metzler matrix~$A$ has a real eigenvalue that is greater than or equal to the real parts of the other eigenvalues of~$A$. This maximum real eigenvalue is denoted by~$\lambda_{\max}(A)$. Let $A\otimes B$ denote the Kronecker product of matrices~$A$ and~$B$. If $A$ and~$B$ are square, then the Kronecker sum of~$A$ and~$B$ is defined by~$A\oplus B = A\otimes I_m + I_n \otimes B$, where $n$ and~$m$ denote the orders of~$A$ and~$B$, respectively. The diagonal matrix having block diagonals~$A_1$, \dots, $A_n$ is denoted by~$\diag(A_1, \dotsc, A_n)$. For a vector $a$ having scalar entries~$a_1$, \dots,~$a_n$, we often use the shorthand notation
\begin{equation*}
D_a = \diag(a_1, \dotsc, a_n).
\end{equation*}

This paper is organized as follows. In Section~\ref{sec:prb}, we formulate the class of optimization problems studied in this paper. Then, in Sections~\ref{sec:h2}--\ref{sec:Hankel}, we present geometric programs for tuning the parameters of positive linear systems constrained by the $H^2$ norm, the $H^\infty$ norm, and the Hankel singular values, respectively. In Section~\ref{sec:robust}, we present a geometric program for tuning the parameters so that the robust stability of the system under structural uncertainties is guaranteed. In Section~\ref{sec:delay}, we show that a class of mixed-constraint parameter tuning problem for delayed positive linear systems reduces to a convex optimization problem. We illustrate the obtained theoretical results in Sections~\ref{sec:simulations2} and~\ref{sec:simulations}. We finally provide the conclusion of the paper as well as some discussions in Section~\ref{sec:conclusion}.

\section{Problem formulation} \label{sec:prb}

In this section, we formulate the problems studied in this paper. Let us consider the linear time-invariant system
\begin{equation*}
\Sigma_\theta:
\begin{cases}
\,\dfrac{dx}{dt} = A(\theta)x + B(\theta)w, 
\\
\,y = C(\theta)x,\end{cases}  
\end{equation*}
which is parametrized by the parameter~$\theta$ belonging to a subset~$\Theta \subset \mathbb R^\nsubtheta$. We suppose that the matrix functions~$A$, $B$, and~$C$ are defined on $\Theta$ and have dimensions $n_x\times n_x$, $n_x\times n_w$, and~$n_y\times n_x$, respectively. 

{To guarantee the (internal) positivity of the system~$\Sigma_{\theta}$, we assume that, for all $\theta \in \Theta$, the matrix~$A(\theta)$ is Metzler and the matrices~$B(\theta)$ and~$C(\theta)$ are nonnegative (see, e.g.,~\cite{Farina2000}). Under these assumptions, for all nonnegative initial condition~$x(0)$ and nonnegative input signal~$u(t)$ ($t\geq 0$), the values of the state $x(t)$ and output $y(t)$ remain nonnegative at every time instant~$t$. Also, we say that the system~$\Sigma_{\theta}$ is internally stable if the matrix~$A(\theta)$ is Hurwitz stable. }

{The parametrized positive model~$\Sigma_{\theta}$ arises in various contexts including drug therapy and leader selection~\cite{Dhingra2017}, as well as dynamical buffer networks~\cite{Rantzer2018} and networked epidemics~\cite{Nowzari2015a,Pastor-Satorras2015a} (see Sections~\ref{sec:simulations} and~\ref{sec:simulations2} for these examples, respectively). }
In this paper, we consider the following general parameter optimization problem:
\begin{equation}\label{eq:generalOptCon}
\begin{aligned}
\minimize_{\theta\in\Theta}\ \ \ & L(\theta)
\\
\subjectto\ \ \,\,& \mbox{$\Sigma_\theta$ is internally stable, }
\\
&J(\Sigma_\theta) \leq \gamma, 
\end{aligned}
\end{equation}
where $\theta$ is the parameter to be tuned, the mapping
\begin{equation*}
L\colon \Theta \to [0, \infty)
\end{equation*}
represents the cost for realizing the parameter~$\theta$, and the constraint~$J(\Sigma_\theta) \leq \gamma$ is our requirement on the system~$\Sigma_\theta$ in terms of a functional~$J$ and a constant~$\gamma$. For example, we allow the functional~$J$ to be the $H^2$ norm of the system~$\Sigma_{\theta}$ defined by 
\begin{equation*}
\norm{\Sigma_\theta}_2 = \sqrt{\int_0^\infty \tr\left(\Phi_\theta(t)\Phi_\theta(t)^\top \right)\,dt}, 
\end{equation*}
where $\Phi_\theta(t) = C(\theta) \exp(A(\theta)t)B(\theta) \in \mathbb{R}^{n_y\times n_w}$ is the impulse response of the system~$\Sigma_\theta$ and~$\tr(\cdot)$ denotes the trace of a matrix. Another functional that we consider is the $\mathcal L^2$-gain ({i.e.}, the $H^\infty$ norm) of the system defined by
\begin{equation*}
\norm{\Sigma_{\theta}}_\infty = \sup_{w \in \mathcal L^2(\mathbb{R}^{n_w})\backslash \{0\}} \frac{\norm{\Phi_\theta * w}_2}{\norm{w}_2}, 
\end{equation*}
where $*$ denotes a convolution product and $\mathcal L^2(\mathbb{R}^n) = \{ f\colon [0, \infty) \to \mathbb{R}^n \mid \int_0^\infty \norm{f(t)}^2\,dt < \infty \}$ denotes the space of Lebesgue-measurable square-integrable functions equipped with the norm $\norm{f}_2 = ( \int_0^\infty \norm{f(t)}^2\,dt )^{1/2}$.

Throughout this paper, we place a certain regularity assumption on the coefficient matrices in the system~$\Sigma_{\theta}$. To state the assumption, we introduce the class of posynomial functions~\cite{Boyd2007}. 

\begin{definition}\label{defN:monoposy}
Let $v_1$, $\dotsc$, $v_n$ denote positive variables and define $v = (v_1, \dotsc, v_n)$. 
\begin{enumerate}
\item We say that a real function $h$ of~$v$ is a \emph{\it monomial} if there exist $c>0$ and~$a_1, \dotsc, a_n \in \mathbb{R}$ such that $h(v) = c v_{\mathstrut 1}^{a_{1}} \dotsm v_{\mathstrut n}^{a_n}$. 

\item We say that a real function~$f$ of~$v$ is a \emph{posynomial} if $f$ is the sum of monomials of~$v$. 
\end{enumerate}
\afterequation
\end{definition}  

Monomials and posynomials are closely related to a class of optimization problems called geometric programs. Given posynomials~$f_0$, \dots, $f_p$ and monomials~$h_1$, \dots, $h_q$, the optimization problem
\begin{equation}\label{eq:scalarGP}
\begin{aligned}
\minimize_{v\in \mathbb{R}^n_{++}}\ \ \ 
&
f_0(v)
\\
\subjectto \ \ \,\,
&
f_i(v)\leq 1,\quad i=1, \dotsc, p, 
\\
&
h_j(v) = 1,\quad j=1, \dotsc, q, 
\end{aligned}
\end{equation}
is called a geometric program~\cite{Boyd2007}. 
It is known~\cite{Boyd2007} that a geometric program can be converted into a convex optimization problem via the logarithmic variable transformation {
\begin{equation}\label{eq:transformation}
v = \exp[z], \ z\in\mathbb{R}^n
\end{equation}
where $\exp[\cdot]$ stands for entrywise exponentiation of a real vector. Specifically, this transformation yields the following equivalent optimization problem
\begin{equation*} 
\begin{aligned}
\minimize_{z\in \mathbb{R}^n}\ \ \ 
&
\log f_0(\exp[z])
\\
\subjectto \ \ \,\,
&
\log f_i(\exp[z])\leq 0,\quad i=1, \dotsc, p, 
\\
&
\log h_j(\exp[z]) = 0,\quad j=1, \dotsc, q, 
\end{aligned}
\end{equation*}
{which can be efficiently solved using, for example, interior-point methods (see~\cite{Boyd2007}, for more details on GP). Specifically, the geometric program \eqref{eq:scalarGP} can be solved with computational cost polynomial in $p$, $q$, and the maximum of the numbers of monomials contained in each of posynomials $f_0$, $\cdots$, $f_p$ \cite[Section~10.4]{Nemirovskii2004}. Furthermore, packages for directly formulating and solving geometric programs are available in various standard softwares including MATLAB, Python, and MOSEK. }

We now state our assumptions on the parametrization of the coefficient matrices in the system~$\Sigma_{\theta}$.  

\begin{assumption}[Coefficient matrices]\label{asm:posy}
The following conditions hold true: 
\begin{enumerate}
\item \label{asm:posy:A} There exists a diagonal matrix function~
\begin{equation}\label{eq:def:R}
R(\theta) = \diag(r_1(\theta),  \dotsc, r_{n_x}(\theta)) 
\end{equation}
having monomial diagonals $r_1(\theta)$, \dots, $r_{n_x}(\theta)$ such that each entry of the matrix
\begin{equation*}
\tilde A(\theta) = A(\theta) + R(\theta)
\end{equation*}
is either a posynomial of~$\theta$ or zero.

\item \label{asm:posy:BC} 
Each entry of the matrices $B(\theta)$ and~$C(\theta)$ is either  a posynomial of~$\theta$ or zero.
\end{enumerate}
\afterequation
\end{assumption}

{\begin{remark}\label{rmk:whatisA}
Assumption~\ref{asm:posy} implicitly limits the parameter set~$\Theta$ to the positive orthant. This limitation allows us to employ the framework of the geometric programming. Also, Assumption~\ref{asm:posy}.\ref{asm:posy:A}) states that the off-diagonals of $A(\theta)$ are either a posynomial or zero, while the diagonals of $A(\theta)$ are \emph{signomials} with at most one negative coefficient (see \cite{Chandrasekaran2016} for the details).
\end{remark}}

Let us also place the following assumptions on the parameter~$\theta$. 

\begin{assumption}[Parameter~$\theta$ and cost~$L(\theta)$]\label{asm:costPosy}
The following conditions hold true: 
\begin{enumerate}
\item \label{ams:item:L} {$L(\theta)$ is a constant shift of a posynomial, {that is, }}there exists a constant $L_0$ such that 
\begin{equation}\label{eq:deftildeL}
\tilde L(\theta) =L(\theta) + L_0
\end{equation}
is a posynomial of~$\theta$. 
\item \label{ams:item:Theta} There exist posynomials $f_1(\theta)$, $\dotsc$, $f_p(\theta)$ such that the constraint set $\Theta$ satisfies
\begin{equation}\label{eq:def:Theta}
\Theta = \{\theta \in \mathbb R^\nsubtheta \mid \theta>0,\,f_1(\theta)\leq 1,\, \dotsc,\, f_p(\theta)\leq 1\}. 
\end{equation}
\afterequation
\end{enumerate}
\end{assumption}

\section{$H^2$ norm-constrained {parameter }optimization}\label{sec:h2}

Let us consider the following $H^2$ norm-constrained {parameter }optimization problem: 
\begin{equation}\label{eq:h2Cont}
\begin{aligned}
\minimize_{\theta\in\Theta}\ \ \  & L(\theta)
\\
\subjectto\ \ \,\, & \mbox{$\Sigma_\theta$ is internally stable, }
\\&\norm{\Sigma_\theta}_2 < \gamma_{\,2}, 
\end{aligned}
\end{equation}
where $\gamma_{\,2} > 0$ is a constant. In this section, we show that this optimization problem can be solved by geometric programming. To state the result, let us introduce the following notations. 
For each $i=1, \dotsc, n_y$ and  $j = 1, \dotsc, n_w$, let $C_i(\theta)$ and~$B_j(\theta)$ denote the $i$th row and~$j$th column of the matrices~$C(\theta)$ and~$B(\theta)$, respectively. Define the $n_x^2$\nobreakdash-dimensional column and row vectors
\begin{equation*}
\begin{aligned}
\tilde B(\theta) &=\sum_{j=1}^{n_w} B_j(\theta)\otimes B_j(\theta), 
\\
 \tilde C(\theta) &= \sum_{i=1}^{n_y} C_i(\theta)\otimes C_i(\theta). 
\end{aligned}
\end{equation*}

{The following theorem states that we can solve the $H^2$ norm-constrained parameter optimization problem if the matrix $R(\theta)$ in Assumption~\ref{asm:posy}.\ref{asm:posy:A}) can be chosen in a specific form.}

\begin{theorem}\label{thm:h2cont}
Assume that there exist a monomial~$r(\theta)$ and a diagonal matrix~$R_0$ with positive diagonals such that {the matrix~$R(\theta)$ given in~\eqref{eq:def:R} satisfies}
\begin{equation}\label{eq:Rfurtherassumption}
R(\theta) = r(\theta)R_0.
\end{equation} 
Then, the solution of the $H^2$ norm-constrained {parameter }optimization problem~\eqref{eq:h2Cont} is given by the solution of the following geometric program:
\begin{subequations}\label{eq:H2GP}
\begin{align}
\minimize_{\mathclap{\theta \in \mathbb R^\nsubtheta_{++},\, \omega \in \mathbb{R}^{n_x^2}_{++}}}\ \ & \tilde L(\theta)\label{eq:confusingPoint}
\\
\subjectto\ \,\,& 
\gamma^{-2}_2 \tilde C(\theta) \omega < 1, \label{eq:h2const1}\\
&\frac{D_{\omega}^{-1}(R_0 \oplus R_0)^{-1}\left[\left(\tilde A(\theta)\oplus \tilde A(\theta)\right)\omega + \tilde B (\theta)\right]}{r(\theta)}< \onev,  
\label{eq:h2const2}\!\!\!\\
&f_i(\theta) \leq 1,\quad i=1, \dotsc, p. 
\label{eq:h2const3}
\end{align}
\end{subequations}
\afterequation
\end{theorem}

\begin{remark}
Geometric programs in standard form do not allow strict inequality constraints appearing in the optimization problem~\eqref{eq:H2GP}. For this reason, in practice, we would relax the strict inequality constraints into non-strict counterparts by, for example, replacing the constraint~\eqref{eq:h2const1} with $\gamma^{-2}_2 \tilde C(\theta) \omega \leq 1-\epsilon$ for a small constant~$\epsilon > 0$. 
\end{remark}

For the proof of Theorem~\ref{thm:h2cont}, we start by showing the following lemma. 

\begin{lemma}[{\cite[Lemma~1]{Briat2012c}}]\label{lem:key}
Let $F \in \mathbb{R}^{n\times n}$, $g \in \mathbb{R}^{n}$, $H  \in \mathbb{R}^{m\times n}$, and~$v  \in \mathbb{R}^m$. Assume that $F$ is Metzler, and~$g$ and~$H $ are nonnegative. The following conditions are equivalent.
\begin{enumerate}
\item $F$ is Hurwitz stable and $-H  F^{-1} g < v$. 
\item There exists a positive vector~$\omega \in \mathbb{R}^n$ such that $H  \mathcal  \omega <v$ and $ F \omega + g < 0 $.
\end{enumerate}
\afterequation
\end{lemma}

We then present the following proposition that characterizes the $H^2$ norm of a positive linear system 
\begin{equation}\label{eq:lti}
\Sigma\colon
\begin{cases}
\,\dfrac{dx}{dt} = Fx + Gw, 
\\
\,y = Hx,\end{cases}  
\end{equation}
where $F$ is a Metzler $n_x\times n_x$ matrix, and~$G$ and~$H$ are  $n_x\times n_w$ and~$n_y\times n_x$ nonnegative matrices. 

\begin{proposition}\label{prop:H2}
Let $\gamma > 0$ be a constant. 
Define the $n_x^2$\nobreakdash-dimensional row and column vectors
 $\tilde H = \sum_{i=1}^{n_y} H_i\otimes H_i$ and~$
\tilde G =\sum_{j=1}^{n_w} G_j\otimes G_j$, 
where $H_i$ and~$G_j$ denote the $i$th row and~$j$th column of the matrices $H$ and~$G$, respectively. Then, the following conditions are equivalent:
\begin{enumerate}
\item $\Sigma$ is internally stable and  $\norm{\Sigma}_2 < \gamma$. 
\item There exists a positive vector~$\omega \in \mathbb{R}^{n_x^2}$ such that 
\begin{equation}\label{eq:H2}
\begin{aligned}
&\tilde H \omega < \gamma^{\,2}, \\
&(F\oplus F) \omega + \tilde G < 0. 
\end{aligned}
\end{equation}
\afterequation
\end{enumerate}
\end{proposition}

\begin{IEEEproof}
Assume that $\Sigma$ is internally stable and $\norm{\Sigma}_2 < \gamma$. By~\cite[Theorem~2]{Ebihara2017a}, we have $\norm{\Sigma}_2^2 = -\tilde H (F\oplus F)^{-1}  \tilde G$. Applying Lemma~\ref{lem:key} to the inequality $-\tilde H (F\oplus F)^{-1}  \tilde G < \gamma$, we can show the existence of a positive vector~$\omega \in \mathbb{R}^{n^2}$ satisfying inequalities in~\eqref{eq:H2}. The other direction of the proof is straightforward, and, therefore, is omitted.
\end{IEEEproof}

Let us prove Theorem~\ref{thm:h2cont}.

\begin{proofof}{Theorem~\ref{thm:h2cont}}
{Assumptions}~\ref{asm:posy} and~\ref{asm:costPosy} show that the optimization problem~\eqref{eq:H2GP} is a geometric program. For example, Assumption~\ref{asm:costPosy} shows that the objective function~$\tilde L(\theta)$ is a posynomial. Also, to confirm that each entry of vector on the left-hand side of the constraint~\eqref{eq:h2const2} is a posynomial, we first notice that any entry of~$D_\omega^{-1}$, $(R_0\oplus R_0)^{-1}$, $\tilde A(\theta)\oplus \tilde A(\theta)$, $\tilde B(\theta)$, and~$1/r(\theta)$ is either a posynomial with the variables~$\theta$ and~$\omega$ or a nonnegative constant. Then, by using the fact that the set of posynomials is closed under addition and multiplications~\cite{Boyd2007}, we can confirm that the constraint~\eqref{eq:h2const2} is indeed written in terms of posynomials.

Let us show that the $H^2$ norm-constrained {parameter }optimization problem~\eqref{eq:h2Cont} reduces to the geometric program~\eqref{eq:H2GP}. Proposition~\ref{prop:H2} implies that the solution of the optimization problem~\eqref{eq:h2Cont} is given by the solution of the following optimization problem: 
\begin{subequations}\label{eq:h2Cont:pre}
\begin{align}
\minimize_{\mathclap{\theta \in \Theta,\, \omega \in \mathbb{R}^{n_x^2}_{++}}}\ \ \ & L(\theta)
\\
\subjectto\ \ \,\,& 
\tilde C(\theta) \omega < \gamma^{\,2}_{\,2}, \label{eq:h2constFirst}\\
&\left(A(\theta)\oplus A(\theta)\right)\omega + \tilde B (\theta)< 0.  \label{eq:h2constSecond}
\end{align}
\end{subequations}
In this optimization problem, the minimization of~$L(\theta)$ is equivalent to minimizing $\tilde L(\theta)$ by the relationship~\eqref{eq:deftildeL}. The constraint~\eqref{eq:h2constFirst} is clearly equivalent to the constraint~\eqref{eq:h2const1}. Furthermore, since we have $A(\theta)\oplus A(\theta) = \tilde A(\theta)\oplus \tilde A(\theta) - r(\theta)(R_0 \oplus R_0)$ and~$D_\omega^{-1}\omega = \onev$, the constraint \eqref{eq:h2constSecond} is equivalent to~\eqref{eq:h2const2}. Finally, \eqref{eq:def:Theta} implies that $\theta \in \Theta$ if and only if constraints~\eqref{eq:h2const3} hold true. Therefore, we conclude that the optimization problem \eqref{eq:h2Cont:pre} reduces to the geometric program~\eqref{eq:H2GP}, as desired.
\end{proofof}

{
\begin{remark}\label{rmk:}
Theorem~\ref{thm:h2cont} has  a few immediate consequences. For example, one can easily confirm that the $H^2$ norm-constrained parameter optimization problem~\eqref{eq:h2Cont} is solvable for all $\gamma_{\,2} \geq \gamma^\star_2$, where $\gamma^\star_2$ is the solution of the following geometric program:
\begin{equation*}
\begin{aligned}
\minimize_{\mathclap{\theta \in \mathbb R^\nsubtheta_{++},\, \omega \in \mathbb{R}^{n_x^2}_{++},\,\gamma_{\,2}>0}}\ \ \ \ \ \ \ \ & \gamma_{\,2}
\\
\subjectto\ \ \ \ \ \ \ \,\,& \mbox{\eqref{eq:h2const1}--\eqref{eq:h2const3}}.
\end{aligned}
\end{equation*}
Similarly, we can show that the following cost-constrained counterpart of the $H^2$ norm-constrained parameter optimization problem~\eqref{eq:h2Cont}:
\begin{equation*}
\begin{aligned}
\minimize_{\theta\in\Theta}\ \ \ & \norm{\Sigma_\theta}_2
\\
\subjectto\ \ \,\, &\mbox{$\Sigma_\theta$ is internally stable, }
\\& L(\theta) \leq \bar L, 
\end{aligned}
\end{equation*}
where $\bar L > 0$ is a given constant, is solved by the following geometric program: 
\begin{equation*}
\begin{aligned}
\minimize_{\mathclap{\theta \in \mathbb R^\nsubtheta_{++},\, \omega \in \mathbb{R}^{n_x^2}_{++},\,\gamma_{\,2}>0}}\ \ \ \ \ \ \ & \gamma_{\,2} 
\\
\subjectto\ \ \ \ \ \ \,\,& \tilde L(\theta)\leq \bar L + L_0,\  \mbox{\eqref{eq:h2const1}--\eqref{eq:h2const3}.}
\end{aligned}
\end{equation*} 
\end{remark}}

\section{$H^\infty$ norm-constrained {parameter }optimization}

In this section, we show that the $H^\infty$ norm-constrained {parameter }optimization problem
\begin{equation}\label{eq:hinfCont}
\begin{aligned}
\minimize_{\theta\in\Theta}\ \ \ & L(\theta)
\\
\subjectto\ \ \,\,& \mbox{$\Sigma_\theta$ is internally stable, }
\\&\norm{\Sigma_\theta}_{\infty} < \gamma_\infty, 
\end{aligned}
\end{equation}
for a positive constant~$\gamma_\infty$ can be solved by  geometric programming, as stated in the following theorem. 

\begin{theorem}\label{thm:Hinf}
The solution of the $H^\infty$ norm-constrained {parameter }optimization problem~\eqref{eq:hinfCont} is given by the solution of the following geometric program:
\begin{subequations}\label{eq:optHinf}
\begin{align}
\minimize_{\mathclap{\substack{\theta\in \mathbb R^\nsubtheta_{++},\\u\in\mathbb{R}^{n_w}_{++},\, v\in\mathbb{R}^{n_y}_{++},\\\xi,\,\zeta \in \mathbb{R}^{n_x}_{++}}}}\ \ \ & \tilde L(\theta)
\\
\subjectto\ \ \,\,&
 \gamma^{-1}_\infty D_v^{-1}C(\theta) \xi  < \onev , \label{eq:hinfconstFirst}
\\
& D_\xi^{-1}R(\theta)^{-1}(\tilde A(\theta)\xi + B(\theta) u) < \onev,
\\
& \gamma^{-1}_\infty D_u^{-1} B(\theta)^\top \zeta < \onev, 
\\
& D_\zeta^{-1} R(\theta)^{-1}(\tilde A(\theta)^\top \zeta + C(\theta)^\top v) < \onev, 
\label{eq:hinfconstsecondLast}
\\
&f_i(\theta) \leq 1, \quad i=1, \dotsc, p \label{eq:hinfconstLast}.
\end{align}
\end{subequations} 
\afterequation
\end{theorem}

{For the proof of Theorem~\ref{thm:Hinf}, we start by recalling the Perron\nobreakdash-Frobenius theorem for Metzler matrices.

\begin{lemma}[{\cite{Horn1990}}]\label{lem:PF}
Let $M$ be an $n\times n$ Metzler matrix and $\gamma$ be a real number. We have $\lambda_{\max}(M) < \gamma$ if and only if there exists a positive vector~$v\in\mathbb{R}^n$ such that $Mv < \gamma\,v$. 
\end{lemma}

Then, we state the following lemma for characterizing the maximum singular value of nonnegative matrices.} 

\begin{lemma}\label{lem:normmat}
Let $M$ be a nonnegative matrix and~$\gamma$ be a positive number. Then, the following conditions are equivalent:
\begin{enumerate}
\item $\norm{{M}} < \gamma$. 
\item There exist positive vectors~$u$ and~$v$ such that  
\begin{subequations}\label{eq:Auvuvuv}
\begin{align}
&{M}u < \gamma\,v,\label{eq:Avsu}
\\
&{M}^\top v < \gamma\,u. \label{eq:Ausv}
\end{align}
\end{subequations}
\end{enumerate}
\afterequation
\end{lemma}

\begin{IEEEproof}
Assume $\norm{M} <\gamma$. This implies $\lambda_{\max}(M^\top M) < \gamma^{\,2}$. By Lemma~\ref{lem:PF}, there exists a positive vector~$u$ such that $M^\top M u < \gamma^{\,2} u$ and, hence, $M^\top M u/\gamma < \gamma\,u$. Therefore, we can take an $\epsilon>0$ such that $M^\top M \frac{u}{\gamma-\epsilon}  < \gamma\,u$. If we define $v =M u/(\gamma-\epsilon)$, then this inequality implies inequality~\eqref{eq:Ausv}. Also, by the definition of the vector~$v$, we have $M u = (\gamma-\epsilon) v < \gamma\,v$, which yields inequality~\eqref{eq:Avsu}.

Conversely, assume that there exist positive vectors~$u$ and~$v$ satisfying \eqref{eq:Auvuvuv}. Then, we have $M^\top M u < \gamma^{\,2} u$. This inequality and Lemma~\ref{lem:PF} show $\lambda_{\max}(M^\top M) < \gamma^{\,2}$.  Hence, we obtain~$\norm{M} < \gamma$, as desired.
\end{IEEEproof}

Using Lemma~\ref{lem:normmat}, we can prove the following proposition for characterizing the $H^\infty$ norm of a positive linear system. 

\begin{proposition}\label{prop:Hinf}
Consider the linear system~$\Sigma$ given in~\eqref{eq:lti}. Let $\gamma > 0$. The following statements are equivalent:
\begin{enumerate}
\item $\Sigma$ is internally stable and~$\norm{\Sigma}_{\infty} < \gamma$. 
\item There exist positive vectors~$u\in\mathbb{R}^{n_w}$, $v \in \mathbb{R}^{n_y}$ and~$\xi, \zeta \in \mathbb{R}^{n_x}$ such that the following inequalities hold true: 
\begin{subequations}
\begin{align}
&H \xi  < \gamma\,v, \label{eq:hinf:analysis:1}
\\
&F\xi + G u < 0,\label{eq:hinf:analysis:2}  
\\
&G^\top \zeta < \gamma\,u,  \label{eq:hinf:analysis:3}
\\
&F^\top \zeta + H^\top v < 0. \label{eq:hinf:analysis:4}
\end{align}
\end{subequations}
\afterequation
\end{enumerate}
\end{proposition}

\begin{IEEEproof}
Assume that $\Sigma$ is internally stable and~$\norm{\Sigma}_{\infty} < \gamma$. Then, by~\cite[Theorem~2]{Tanaka2013}, we have $\norm{\hat M(0)} < \gamma$ for the transfer function~$\hat M(s) = H(sI-F)^{-1}G$. Since $\hat M(0) = -HF^{-1}G$, Lemma~\ref{lem:normmat} shows the existence of positive vectors~$u$ and~$v$  such that
\begin{align}
&{- HF^{-1}G u} < \gamma\,v, \label{eq:hinf:inequalities1}
\\
&{- G^\top (F^{\top})^{-1} H^\top v}  < \gamma\,u.   \label{eq:hinf:inequalities2}
\end{align}
Since $F$ is Hurwitz stable, we can apply Lemma~\ref{lem:key} to inequality~\eqref{eq:hinf:inequalities1} to show the existence of a positive vector~$\xi$ for which inequalities~\eqref{eq:hinf:analysis:1} and~\eqref{eq:hinf:analysis:2} hold true. Similarly, applying Lemma~\ref{lem:key} to inequality~\eqref{eq:hinf:inequalities2}, we can show the existence of a positive vector~$\zeta$ satisfying inequalities~\eqref{eq:hinf:analysis:3} and~\eqref{eq:hinf:analysis:4}. The proof of the other direction is omitted.
\end{IEEEproof}

We can now prove Theorem~\ref{thm:Hinf}.

\begin{proofof}{Theorem~\ref{thm:Hinf}}
Proposition~\ref{prop:Hinf} implies that the solution of the $H^\infty$ norm-constrained optimization problem~\eqref{eq:hinfCont} is given by the solution of the following optimization problem: 
\begin{equation*}
\begin{aligned}
\minimize_{\mathclap{\substack{\theta\in \Theta,\,u\in\mathbb{R}^{n_w}_{++},\, v\in\mathbb{R}^{n_y}_{++},\\\xi,\,\zeta \in \mathbb{R}^{n_x}_{++}}}}\ \ \ \ \ \ \ & L(\theta)
\\
\subjectto\ \ \ \ \ \ \,\,&
 C(\theta) \xi  < \gamma_\infty v, 
\\
& A(\theta)\xi + B(\theta) u <0,
\\
& B(\theta)^\top \zeta < \gamma_\infty u, 
\\
& A(\theta)^\top \zeta + C(\theta)^\top v < 0.
\end{aligned}
\end{equation*}
An algebraic manipulation and equalities~\eqref{eq:deftildeL} and~\eqref{eq:def:Theta} show that this optimization problem is equivalent to the optimization problem~\eqref{eq:optHinf}. Furthermore, Assumptions~\ref{asm:posy} and~\ref{asm:costPosy} show that the optimization problem~\eqref{eq:optHinf} is indeed a geometric program. This completes the proof of the theorem.
\end{proofof}

{A few remarks are in order. First, as stated in {Remark~\ref{rmk:} } for the case of the $H^2$ norm, we can derive geometric programs for 1) finding the minimum achievable $H^\infty$ norm of the system and 2) solving a cost-constrained $H^\infty$ norm optimization problem. Since their derivations are straightforward, we do not explicitly state them in this paper. We also remark that, }by using Theorem~\ref{thm:Hinf} as well as Theorem~\ref{thm:h2cont}, we can show that a class of mixed $H^2$/$H^\infty$ optimization problems for positive linear systems reduces to a geometric program. Let us consider the following optimization problem:
\begin{subequations}\label{eq:hinfh2budget}
\begin{align}
\minimize_{\theta\in\Theta}\ \ \ & L(\theta)
\\
\subjectto\ \ \,\, &\mbox{$\Sigma_\theta$ is internally stable, } \label{eq:hinfh2budget:0}
\\
& \alpha(\norm{\Sigma_\theta}_2, \norm{\Sigma_\theta}_\infty) < \gamma, \label{eq:hinfh2budget:1}
\end{align}
\end{subequations}
where $\alpha\colon \mathbb{R}_{++}^2\to \mathbb{R}_{++}$ is a function representing a trade-off between the $H^2$ and~$H^\infty$ norm of the system. Let us place the following assumption on the trade-off function. 

\begin{assumption}\label{asm:alpha}
The function $\alpha$ is a posynomial, and nondecreasing with respect to each variable. 
\end{assumption}

Examples of the function~{$\alpha(\norm{\Sigma_\theta}_2, \norm{\Sigma_\theta}_\infty)$ }satisfying these assumptions include the sum $\norm{\Sigma_\theta}_2 + \norm{\Sigma_\theta}_\infty$ and the product~$\norm{\Sigma_\theta}_2\norm{\Sigma_\theta}_\infty$. Under this assumption, the following corollary shows that the solution of the mixed $H^2$/$H^\infty$ optimization problem~\eqref{eq:hinfh2budget} is obtained by geometric programming.

\begin{corollary}\label{cor:mixedH2Hinfty}
If there exist a monomial $r(\theta)$ and a diagonal matrix~$R_0$ with positive diagonals satisfying~\eqref{eq:Rfurtherassumption}, then the solution of the mixed $H^2$/$H^\infty$ norm-constrained {parameter }optimization problem~\eqref{eq:hinfh2budget} is given by the solution of the following geometric program:
\begin{subequations}\label{eq:optmix}
\begin{align}
\minimize_{\mathclap{\substack{\theta \in \mathbb R^\nsubtheta_{++},\, \omega \in \mathbb{R}^{n_x^2}_{++},\\u\in\mathbb{R}^{n_w}_{++},\, v\in\mathbb{R}^{n_y}_{++},\,\xi,\,\zeta \in \mathbb{R}^{n_x}_{++}\\\gamma_{\,2},\,\gamma_\infty>0}}}\ \ \ \ \ \ \ & \tilde L(\theta)
\\
\subjectto\ \ \ \ \ \ \,\,& \gamma^{-1}\alpha(\gamma_{\,2}, \gamma_\infty) < 1, \label{eq:optmix:const1}
\\
&\mbox{\eqref{eq:h2const1}, \eqref{eq:h2const2}, \eqref{eq:hinfconstFirst}--\eqref{eq:hinfconstsecondLast}},  \label{eq:optmix:const2}
\\
&f_i(\theta) \leq 1,\quad i=1, \dotsc, p. 
\end{align}
\end{subequations}
\afterequation
\end{corollary}

\newcounter{MYtempeqncnt}
\begin{figure*}[!t]
\normalsize
\setcounter{MYtempeqncnt}{\value{equation}}
\setcounter{equation}{23}
\begin{subequations}\label{eq:HankelGP}
\begin{align}
\minimize_{\mathclap{\substack{\theta \in \mathbb R^\nsubtheta_{++},\,v\in\mathbb{R}^{n_x}_{++}\\\, \omega_1 \in \mathbb{R}^{n_x^2n_w}_{++},\, \omega_2 \in \mathbb{R}^{n_x^2n_y}_{++}} } }\ \ \ \ \ \ \ & \tilde L(\theta)
\\
\subjectto\ \ \ \ \ \ \,\,& 
\gamma^{-2} D_v^{-1} \bar B_1(\theta) \omega_1 < \onev, \label{eq:Hankelconst1}\\
&D_{\omega_1}^{-1}(R_0 \oplus O_{n_w} \oplus R_0)^{-1}\frac{\bar B_2(\theta) \bar C_1(\theta) \omega_2 + (\tilde A(\theta) \oplus  O_{n_w} \oplus  \tilde A(\theta)^\top)\omega_1}{r(\theta)} < \onev,  
\label{eq:Hankelconst2}\\
&
D_{\omega_2}^{-1}(R_0 \oplus O_{n_y} \oplus R_0)^{-1}\frac{(\tilde A(\theta)^\top \oplus  O_{n_y} \oplus  \tilde A(\theta))\omega_2+\bar B_2(\theta)v}{r(\theta)} < \onev,  
\label{eq:Hankelconst3}\\
&f_i(\theta) \leq 1,\quad i=1, \dotsc, p. 
\label{eq:Hankelconst4}
\end{align}
\end{subequations}
\setcounter{equation}{\theMYtempeqncnt}
\hrulefill
\vspace*{4pt}
\end{figure*}

\begin{IEEEproof}
The optimization problem~\eqref{eq:optmix} is a geometric program by Assumptions~\ref{asm:posy}, \ref{asm:costPosy}, and~\ref{asm:alpha}. Let $\theta \in \Theta$ and~$\gamma > 0$ be arbitrary. We need to show that the constraints~\eqref{eq:hinfh2budget:0} and~\eqref{eq:hinfh2budget:1} hold true if and only if there exist vectors~$\omega \in \mathbb{R}^{n_x^2}_{++}$, $u\in\mathbb{R}^{n_w}_{++}$, $v\in\mathbb{R}^{n_y}_{++}$, and~$\xi,\zeta \in \mathbb{R}^{n_x}_{++}$ as well as positive constants~$\gamma_{\,2}$ and~$\gamma_\infty$ satisfying constraints~\eqref{eq:optmix:const1} and~\eqref{eq:optmix:const2}.

Assume that \eqref{eq:hinfh2budget:0} and \eqref{eq:hinfh2budget:1} hold true. Then, by the continuity of posynomials, there exist constants~$\gamma_{\,2}$ and~$\gamma_\infty$ satisfying $\norm{\Sigma_\theta}_2 < \gamma_{\,2}$, $\norm{\Sigma_\theta}_\infty < \gamma_\infty$, and \eqref{eq:optmix:const1}. Then, Propositions~\ref{prop:H2} and~\ref{prop:Hinf} show the existence of the vectors~$\omega \in \mathbb{R}^{n_x^2}_{++}$, $u\in\mathbb{R}^{n_w}_{++}$, $v\in\mathbb{R}^{n_y}_{++}$, and~$\xi,\zeta \in \mathbb{R}^{n_x}_{++}$ satisfying \eqref{eq:optmix:const2} as well. Conversely, assume that there exist $\omega \in \mathbb{R}^{n_x^2}_{++}$, $u\in\mathbb{R}^{n_w}_{++}$, $v\in\mathbb{R}^{n_y}_{++}$, and~$\xi,\zeta \in \mathbb{R}^{n_x}_{++}$ as well as positive constants~$\gamma_{\,2}$ and~$\gamma_\infty$ satisfying \eqref{eq:optmix:const1} and~\eqref{eq:optmix:const2}. Then, Propositions~\ref{prop:H2} and~\ref{prop:Hinf} show that the system~$\Sigma_{\theta}$ is internally stable and satisfies $\norm{\Sigma_\theta}_2 < \gamma_{\,2}$ and~$\norm{\Sigma_\theta}_\infty < \gamma_\infty$. Furthermore, since $\alpha$ is non-decreasing with respect to both arguments, we obtain~$\alpha(\norm{\Sigma_\theta}_2, \norm{\Sigma_\theta}_\infty) \leq \alpha(\gamma_{\,2}, \gamma_\infty) < \gamma$ from \eqref{eq:optmix:const1}, as desired. This completes the proof of the corollary.
\end{IEEEproof}

\section{Hankel singular values-constrained {parameter }optimizations}\label{sec:Hankel}

In this section, we show that the {parameter }optimization problem~\eqref{eq:generalOptCon} reduces to a geometric program when constrained by system norms induced from Hankel singular values. Assume that the system~$\Sigma_\theta$ is internally stable. The Hankel singular values of~$\Sigma_{\theta}$,  denoted by $$\sigma_1(\theta)\geq \cdots \geq \sigma_{n_x}(\theta)\geq 0,$$ are defined as the singular values of the Hankel operator associated with the system~$\Sigma_{\theta}$ (see, e.g.,~\cite{Glover1984a}). It is well known that $\sigma_i(\theta) = \sqrt{\lambda_{i}(W_O(\theta)W_C(\theta))}$ holds for all $i=1, \dotsc, n_x$, where $W_C(\theta)$ and $W_O(\theta)$ denote the controllability and observability Grammians defined by
\begin{equation*}
\begin{aligned}
W_C(\theta) &= \int^{\infty}_{0} e^{A(\theta)t} B(\theta) B^\top(\theta) e^{A(\theta)^\top t}\,dt, 
\\
W_O(\theta) &= \int^{\infty}_{0} e^{A^\top(\theta)t} C(\theta)^\top C(\theta) e^{A(\theta) t}\,dt, 
\end{aligned}
\end{equation*}
and $\lambda_1(W_O(\theta)W_C(\theta)) \geq \cdots \geq \lambda_{n_x}(W_O(\theta)W_C(\theta)) \geq 0$ denote the eigenvalues of the matrix~$W_O(\theta)W_C(\theta)$.

{The Hankel singular values induce several interesting system norms. An important example is the Hankel norm
\begin{equation*}
\norm{\Sigma_{\theta}}_{\mathcal H} = \sigma_1(\theta). 
\end{equation*}
Another example is the Schatten $p$-norm (see, e.g., \cite{Opmeer2015}) defined by 
\begin{equation*}
\norm{\Sigma_{\theta}}_{S_p} = \biggl(\sum_{i=1}^{n_x} \sigma_i(\theta)^p\biggr)^{1/p}
\end{equation*}
for a positive integer~$p$, which generalizes the Hilbert-Schmidt norm~$\sqrt{\sum_{i=1}^{n_x} \sigma_i(\theta)^2}$ and the nuclear norm~$\sum_{i=1}^{n_x} \sigma_i(\theta)$.}

In this section, we first show that the Hankel norm-constrained {parameter }optimization problem
\begin{equation}\label{eq:HankelCont}
\begin{aligned}
\minimize_{\theta\in\Theta}\ \ \  & L(\theta)
\\
\subjectto\ \ \,\, & \mbox{$\Sigma_\theta$ is internally stable, }
\\&\norm{\Sigma_\theta}_{\mathcal H} < \gamma, 
\end{aligned}
\end{equation}
can be solved by geometric programming. To state the result, we define the matrix functions
\begin{equation}\label{eq:def:checks}
\begin{aligned}
\check B_{1}(\theta) &= \sum_{k=1}^{n_x}
\bigl(e_k^\top B(\theta)\bigr)\otimes e_k^\top,
\\
\check B_{2}(\theta) &= \sum_{k=1}^{n_x}
e_k \otimes \bigl(B(\theta)^\top e_k\bigr),
\\
\check C_{1}(\theta) &= \sum_{k=1}^{n_x}
\bigr(e_k^\top C(\theta)^\top\bigl)\otimes e_k^\top,
\\
\check C_{2}(\theta) &= \sum_{k=1}^{n_x}
e_k \otimes \bigl(C(\theta) e_k\bigr), 
\end{aligned}
\end{equation}
where $\{e_1, \dotsc, e_{n_x}\}$ is the canonical basis of~$\mathbb{R}^{n_x}$.
Then, let us define
\begin{equation}\label{eq:def:bars}
\begin{aligned}
\bar B_1(\theta) &= \begin{bmatrix}
e_1^\top \otimes \check B_{1}(\theta)
\\
\vdots 
\\
e_{n_x}^\top \otimes \check B_{1}(\theta)
\end{bmatrix},
\\
\bar B_2(\theta) &= \begin{bmatrix}
\check B_{2}(\theta) \otimes e_1 & \cdots & \check B_{2}(\theta) \otimes e_{n_x}
\end{bmatrix}, 
\\
\bar C_1(\theta) &= 
\begin{bmatrix}
e_1^\top \otimes \check C_{1}(\theta)
\\
\vdots 
\\
e_{n_x}^\top \otimes \check C_{1}(\theta)
\end{bmatrix}, 
\\
\bar C_2(\theta) &= \begin{bmatrix}
e_1^\top \otimes \check C_{2}(\theta)
&
\cdots
&
e_{n_x}^\top \otimes \check C_{2}(\theta)
\end{bmatrix}.
\end{aligned}
\end{equation}

\begin{theorem}\label{thm:Hankel}
Assume that there exist a monomial $r(\theta)$ and a diagonal matrix~$R_0$ with positive diagonals such that the matrix~$R(\theta)$ given in~\eqref{eq:def:R} satisfies \eqref{eq:Rfurtherassumption}. 
Then, the solution of the Hankel norm-constrained {parameter }optimization problem~\eqref{eq:HankelCont} is given by the solution of the geometric program~\eqref{eq:HankelGP}.
\stepcounter{equation}
\end{theorem}

For the proof of Theorem~\ref{thm:Hankel}, we state the following extension of Lemma~\ref{lem:key}. 

\begin{lemma}\label{lem:keykey}
Let $q$ be an even integer. For each $i=1, \dotsc, q$, let $F_i \in \mathbb{R}^{n_{i}\times n_i}$ and 
$H_i \in \mathbb{R}^{{n_{i-1}}\times n_i}$ be real matrices. Let $v  \in \mathbb{R}^{n_0}$ and $g \in \mathbb{R}^{n_q}$ be real vectors. Assume that $F_1$, \dots,~$F_q$ are Metzler and $H_1$, \dots, $H_q$, $g$ are nonnegative. The following conditions are equivalent.
\begin{enumerate}
\item The matrices $F_1$, \dots,~$F_q$ are Hurwitz stable and 
\begin{equation}\label{eq:CABCAB}
(H_1F_1)\dotsm (H_qF_q)  g < v. 
\end{equation}

\item There exist positive vectors~$\omega_i \in \mathbb{R}^{n_i}$ ($i=1, \dotsc, q$) such that the following system of inequalities hold true: 
\begin{equation}\label{eq:HF}
\begin{aligned}
H_1 \omega_1 &<v, 
\\
F_i\omega_i + H_{i+1} \omega_{i+1} &< 0,\ \mbox{($i=1, \dotsc, q-1$)}
\\
F_q \omega_q + g &< 0.
\end{aligned}
\end{equation}
\end{enumerate}
\afterequation
\end{lemma}

\begin{IEEEproof}
If inequality~\eqref{eq:CABCAB} holds true, then 
applying Lemma~\ref{lem:key} to \eqref{eq:CABCAB} iteratively $q$ times show the existence of positive vectors~$\omega_i \in \mathbb{R}^{n_i}$ satisfying the inequalities in~\eqref{eq:HF}. The proof of the opposite direction is straightforward and, therefore, is omitted. 
\end{IEEEproof}

We also collect basic facts on Kronecker products and sums in the following lemma.

\begin{lemma}[{\cite{Brewer1978}}]\label{lem:kron}
The following claims hold true.
\begin{enumerate}
\item \label{item:lem:kronsumprod}Let $M$ be a real square matrix. Then, we have $\exp(M) \otimes \exp(M) = \exp({M}\oplus {M})$. 
\item \label{item:lem:kronsumspectrum} Let $M$ and $N$ be $n\times n$ real square matrices having eigenvalues $\{\mu_i \}_{i=1}^n$ and $\{\sigma_i \}_{i=1}^n$, respectively. Then, the set of the eigenvalues of $M\oplus N$ coincides with $\{\mu_i + \sigma_j \}_{i, j=1}^n$. 
\item \label{item:lem:prodDist}Let $M_1$, $M_2$, $N_1$, and $N_2$ be matrices. Assume that the products $M_1M_2$ and $N_1N_2$ are well-defined. Then, $(M_1M_2)\otimes (N_1 N_2) = (M_1\otimes N_1)(M_2\otimes N_2)$. 
\end{enumerate}
\end{lemma}

Let us prove Theorem~\ref{thm:Hankel}. 

\begin{figure*}[!t]
\normalsize
\setcounter{MYtempeqncnt}{\value{equation}}
\setcounter{equation}{30}
\begin{subequations}\label{eq:SchattenGP}
\begin{align}
\minimize_{\mathclap{\substack{\theta \in \mathbb R^\nsubtheta_{++},\,\gamma_i > 0,\\\omega_{i, 2k-1} \in \mathbb{R}^{n_x^2n_w}_{++},\\\omega_{i, 2k} \in \mathbb{R}^{n_x^2n_y}_{++}} } }\ \ \ \ \ \ \ & \tilde L(\theta)
\\
\subjectto\ \ \ \ \ \ \,\,& 
\gamma^{-p}\sum_{i=1}^{n_x}\gamma_i < 1, \label{eq:Schatten0}
\\
& \gamma^{-1}_i e_i^\top \bar C_1(\theta) \omega_1 < 1, \label{eq:Schatten1}\\
&D_{\omega_{i, 2k-1}}^{-1}(R_0 \oplus O_{n_w} \oplus R_0)^{-1}\frac{\bar B_2(\theta) \bar C_1(\theta) \omega_{i, 2k} + (\tilde A(\theta) \oplus  O_{n_w} \oplus  \tilde A(\theta)^\top)\omega_{i, 2k-1}}{r(\theta)} < \onev,\quad \mbox{$k=1, 2, \dotsc, \frac{p}{2}$}
\label{eq:Schatten2}\\
&D_{\omega_{i, 2k}}^{-1}(R_0 \oplus O_{n_y} \oplus R_0)^{-1}\frac{\bar C_2(\theta) \bar B_1(\theta) \omega_{i, 2k+1} + (\tilde A(\theta)^\top \oplus  O_{n_y} \oplus  \tilde A(\theta))\omega_{i, 2k}}{r(\theta)} < \onev,\quad \mbox{$k=1, 2, \dotsc, \frac{p}{2}-1$}
\label{eq:Schatten3}\\
&
D_{\omega_{i, p}}^{-1}(R_0 \oplus O_{n_y} \oplus R_0)^{-1}\frac{(\tilde A(\theta)^\top \oplus  O_{n_y} \oplus  \tilde A(\theta))\omega_{i, p}+\bar B_2(\theta)e_i}{r(\theta)} < \onev,  
\label{eq:Schatten4}\\
&f_i(\theta) \leq 1,\quad i=1, \dotsc, p. 
\label{eq:Schatten5}
\end{align}
\end{subequations}
\setcounter{equation}{\theMYtempeqncnt}
\hrulefill
\vspace*{4pt}
\end{figure*}

\begin{IEEEproof}[Proof of Theorem~\ref{thm:Hankel}] 
Notice that the matrix functions~$\bar B_1$, $\bar B_2$, $\bar C_1$, and $\bar C_2$ are posynomials with the variable~$\theta$ by equations~\eqref{eq:def:checks} and~\eqref{eq:def:bars}. Therefore, it is easy to see that the constraints~\eqref{eq:Hankelconst1}--\eqref{eq:Hankelconst4} are in terms of posynomials under Assumptions~\ref{asm:posy} and~\ref{asm:costPosy}. Therefore, the optimization problem~\eqref{eq:HankelGP} is indeed a geometric program. Hence, to prove Theorem~\ref{thm:Hankel}, it is sufficient to show that $\Sigma_{\theta}$ is internally stable and satisfies $\norm{\Sigma_{\theta}}_{\mathcal H}<\gamma$ if and only if there exist positive vectors~$v\in\mathbb{R}^{n_x}$, $\omega_1 \in \mathbb{R}^{n_x^2n_w}$, and~$\omega_2 \in \mathbb{R}^{n_x^2n_y}$ such that inequalities \eqref{eq:Hankelconst1}--\eqref{eq:Hankelconst3} hold true.

Assume that $\Sigma_{\theta}$ is internally stable and satisfies $\norm{\Sigma_{\theta}}_{\mathcal H}<\gamma$. Let us derive alternative expressions for the Grammian matrices. Let $i, j, k, \ell \in \{1, \dotsc, n_x\}$ and define the scalar function
\begin{equation*}
w_{C, ijkl}(\theta) = e_i^\top e^{A(\theta)t} e_k(e_k^\top B(\theta)) (B(\theta)^\top e_\ell) e_\ell^\top e^{A(\theta)^\top t} e_j. 
\end{equation*}
Since 
$B(\theta)B(\theta)^\top = \sum_{k, \ell = 1}^{n_x} e_k(e_k^\top B(\theta)) (B(\theta)^\top e_\ell) e_\ell^\top$, we can write the $(i, j)$th entry of the controllability Grammian $W_C(\theta)$ as 
\begin{equation}\label{eq:WCij}
[W_C(\theta)]_{i, j } = \sum_{k, \ell=1}^{n_x} \int_0^\infty w_{C, ijkl}(\theta) \,dt. 
\end{equation}
Since the scalar $(e_k^\top B(\theta)) (B(\theta)^\top e_\ell)$ equals its transpose $(e_\ell^\top B(\theta)) (B(\theta)^\top e_k)$, we can rewrite the function~$w_{C, ijkl}(\theta)$ as 
\begin{equation*}
\begin{aligned}
&\,w_{C, ijkl}(\theta) 
\\
=&\, \bigl[e_i^\top e^{A(\theta)t} e_k\bigr] \bigl[(e_\ell^\top B(\theta)) (B(\theta)^\top e_k)\bigr] \bigl[ e_\ell^\top e^{A(\theta)^\top t} e_j\bigr]
\\
=&\,
\bigl[e_i^\top e^{A(\theta)t} e_k\bigr]\otimes \bigl[(e_\ell^\top B(\theta)) I_{n_w} (B(\theta)^\top e_k)\bigr] \otimes \bigl[ e_\ell^\top e^{A(\theta)^\top t} e_j\bigr]
\end{aligned}
\end{equation*}
by using the fact that the product of scalars equals the Kronecker product of the scalars. Then, Lemma~\ref{lem:kron}.\ref{item:lem:prodDist} shows that
\begin{equation}\label{eq:star}
\begin{multlined}[.875\linewidth]
w_{C, ijkl}(\theta) =
\bigl[e_i^\top\otimes(e_\ell^\top B(\theta)) \otimes e_\ell^\top\bigr]
\\
\bigl[e^{A(\theta)t}\otimes I_{n_w}\otimes e^{A(\theta)^\top t}\bigr]
\bigl[e_k \otimes (B(\theta)^\top e_k)\otimes e_j\bigr]. 
\end{multlined}
\end{equation}
We then use Lemma~\ref{lem:kron}.\ref{item:lem:kronsumprod} twice to obtain
\begin{equation}\label{eq:eAIA}
\begin{aligned}
e^{A(\theta)t}\otimes I_{n_w}\otimes e^{A(\theta)^\top t} 
&= e^{A(\theta)t}\otimes e^{O_{n_w}t}\otimes e^{A(\theta)^\top t} 
\\
&= e^{(A(\theta)\oplus O_{n_w}\oplus A(\theta)^\top)t}.
\end{aligned}
\end{equation} 
Since the matrix~$A(\theta)$ is Hurwitz stable by our assumption, the eigenvalues of the Kronecker sum~$A(\theta)\oplus O_{n_w}\oplus A(\theta)^\top$ have negative real part by Lemma~\ref{lem:kron}.\ref{item:lem:kronsumspectrum}. Therefore, from \eqref{eq:eAIA} we obtain~$\int_0^\infty e^{A(\theta)t}\otimes I_{n_w}\otimes e^{A(\theta)^\top t} \,{dt} = -(A(\theta)\oplus O_{n_w}\oplus A(\theta)^\top)^{-1}$. Hence, equations~\eqref{eq:WCij} and~\eqref{eq:star} show that
\begin{equation*}
\begin{aligned}
&\,[W_C(\theta)]_{i, j }
\\
=&\, {-\sum_{k, \ell=1}^{n_x}
\bigl(e_i^\top\otimes(e_\ell^\top B(\theta)) \otimes e_\ell^\top\bigr)
\bigl(A(\theta)\oplus O_{n_w}\oplus A(\theta)^\top\bigr)^{-1}}
\\
&\hspace{5cm} 
 \bigl(e_k \otimes (B(\theta)^\top e_k)\otimes e_j\bigr)
\\
=&\,
-\bigl(e_i^\top\otimes \check B_1(\theta)\bigr)
\bigl(A(\theta)\oplus O_{n_w}\oplus A(\theta)^\top\bigr)^{-1} 
\bigl(\check B_2(\theta)\otimes e_j\bigr), 
\end{aligned}
\end{equation*}
which yields $W_C(\theta) = -\bar B_1(\theta) (A(\theta)\oplus O_{n_w}\oplus A(\theta)^\top)^{-1} \bar B_2(\theta)$. Similarly, we can show that the observability Grammian admits the representation~$W_O(\theta) = -\bar C_1(\theta) (A(\theta)^\top \oplus O_{n_y}\oplus A(\theta))^{-1}  \bar C_2(\theta)$. Now, since $\lambda_{\max}(W_O(\theta) W_C(\theta)) = \lambda_{1}(W_O(\theta) W_C(\theta)) =\norm{\Sigma_{\theta}}_{\mathcal H}^2<\gamma^{\,2}$,  Lemma~\ref{lem:PF} shows the existence of a positive vector $v\in\mathbb{R}^{n_x}$ such that
\begin{equation*}
\begin{multlined}[.8\linewidth]
\bar C_1(\theta) \bigl(A(\theta)^\top \oplus O_{n_y}\oplus A(\theta)\bigr)^{-1}  \bar C_2(\theta) \\\bar B_1(\theta) \bigl(A(\theta)\oplus O_{n_w}\oplus A(\theta)^\top\bigr)^{-1}  \bar B_2(\theta) v < \gamma^{\,2} v.  
\end{multlined}
\end{equation*}
Hence, Lemma~\ref{lem:keykey} shows the existence of positive vectors~$\omega_1 \in \mathbb{R}^{n_x^2n_w}$ and~$\omega_2 \in \mathbb{R}^{n_x^2n_y}$ such that
\begin{equation*}
\begin{aligned}
\bar C_1(\theta) \omega_1 &<\gamma^{\,2} v , 
\\
(A(\theta)^\top \oplus O_{n_y}\oplus A(\theta)) \omega_1 + \bar C_2(\theta) \bar B_1(\theta)  \omega_2 &< 0, 
\\
(A(\theta)\oplus O_{n_w}\oplus A(\theta)^\top) \omega_2 + \bar B_2(\theta) v &< 0. 
\end{aligned}
\end{equation*}
Finally, an algebraic manipulation shows that these inequalities are equivalent to the constraints~\eqref{eq:Hankelconst1}--\eqref{eq:Hankelconst3}, as desired. 

We can similarly prove that the existence of positive vectors~$v\in\mathbb{R}^{n_x}$, $\omega_1 \in \mathbb{R}^{n_x^2n_w}$, and~$\omega_2 \in \mathbb{R}^{n_x^2n_y}$ satisfying \eqref{eq:Hankelconst1}--\eqref{eq:Hankelconst3} shows the internal stability of~$\Sigma_{\theta}$ and inequality~$\norm{\Sigma_{\theta}}_{\mathcal H}<\gamma$. The details are omitted.
\end{IEEEproof}

Let us then consider the following Schatten $p$ norm-constrained {parameter }optimization problem:
\begin{equation}\label{eq:SchattenpOptimization}
\begin{aligned}
\minimize_{\theta\in\Theta}\ \ \  & L(\theta)
\\
\subjectto\ \ \,\, & \mbox{$\Sigma_\theta$ is internally stable, }
\\&\norm{\Sigma_\theta}_{S_p} < \gamma, 
\end{aligned}
\end{equation}
for a constant~$\gamma> 0$. The following theorem shows that this optimization problem can be solved by geometric programming under the assumption that $p$ is an even integer, which covers the interesting case of the Hilbert-Schmidt norm.

\begin{theorem}\label{thm:Schatten}
Suppose that $p$ is an even integer. Assume that there exist a monomial $r(\theta)$ and a diagonal matrix~$R_0$ with positive diagonals such that the matrix~$R(\theta)$ given in~\eqref{eq:def:R} satisfies \eqref{eq:Rfurtherassumption}. Then, the solution of the Schatten $p$ norm-constrained {parameter }optimization problem~\eqref{eq:SchattenpOptimization} is given by the solution of the geometric program~\eqref{eq:SchattenGP}.
\stepcounter{equation} 
\end{theorem}

\begin{IEEEproof}
Suppose that $\Sigma_{\theta}$ is internally stable. Let us first show that $\norm{\Sigma_\theta}_{S_p} < \gamma$ if and only if there exist positive numbers $\gamma_1$, \dots, $\gamma_{n_x}$ satisfying \eqref{eq:Schatten0} and 
\begin{equation}\label{eq:eWWi<gamma}
e_i^\top (W_O(\theta) W_C(\theta))^{p/2} e_i < \gamma_i
\end{equation}
for all $i=1, \dotsc, n_x$. Assume $\norm{\Sigma_\theta}_{S_p} < \gamma$. Since the definition of the Schatten $p$-norm shows
\begin{equation}\label{eq:SchattenConseq}
\norm{\Sigma_\theta}_{S_p} = \Bigl[\tr \bigl((W_O(\theta) W_C(\theta))^{p/2}\bigr)\Bigr]^{1/p}, 
\end{equation}
we obtain~$\tr ((W_O(\theta) W_C(\theta))^{p/2}) < \gamma^p$. From this inequality, we can take positive numbers $\gamma_1$, \dots, $\gamma_{n_x}$ such that $[(W_O(\theta) W_C(\theta))^{p/2}]_{ii} < \gamma_i$ for all $i$ and $\gamma_1 + \cdots + \gamma_{n_x} < \gamma^p$, as desired. On the other hand, if there exist positive numbers $\gamma_1$, \dots, $\gamma_{n_x}$ such that \eqref{eq:Schatten0} and~\eqref{eq:eWWi<gamma} hold true, then \eqref{eq:SchattenConseq} shows $\norm{\Sigma_\theta}_{S_p} = (\sum_{i=1}^n [(W_O(\theta) W_C(\theta))^{p/2}]_{ii})^{1/p} < (\sum_{i=1}^{n_x}\gamma_i)^{1/p} < \gamma$, as desired.

From the above observation, to prove the theorem, we need to show that inequality~\eqref{eq:eWWi<gamma} holds true if and only if there exist positive vectors~$\omega_{i, 2k-1} \in \mathbb{R}^{n_x^2n_w}$ and~$\omega_{i, 2k} \in \mathbb{R}^{n_x^2n_y}$ ($k=1, \dotsc, p/2$) satisfying constraints~\eqref{eq:Schatten1}--\eqref{eq:Schatten4}. We can show this equivalence by applying Lemma~\ref{lem:keykey} to the inequality~\eqref{eq:eWWi<gamma} because the product on the left hand side of~\eqref{eq:eWWi<gamma} is rewritten as $e_i^\top (W_O(\theta) W_C(\theta))^{p/2}e_i = (H_1F_1)\dotsm (H_pF_p)  g$ for the matrices~$H_i, F_i$ given by 
\begin{equation*}
\begin{aligned}
& H_1 = e_i^\top \bar C_1(\theta),
\\
&H_3 = H_5 = \cdots = H_{p-1} = \bar B_2(\theta) \bar C_1(\theta), \\
& H_2 = H_4 = \cdots = H_p = \bar C_2(\theta)\bar B_1(\theta), 
\\
& F_1 = F_3 = \cdots = F_{p-1} = A(\theta)^\top \oplus O_{n_y} \oplus A(\theta), 
\\
& F_2 = F_4 = \cdots = F_p = A(\theta) \oplus O_{n_w} \oplus A(\theta)^\top, 
\end{aligned}
\end{equation*}
and the vector~$g = \bar B_2(\theta) e_i$. The further details of the proof is omitted.
\end{IEEEproof}

\section{Stabilization under structured uncertainty}\label{sec:robust}

In this section, we show that a class of robust stabilization problems under structural uncertainties can be solved by geometric programming. Throughout this section, we place the following assumption for simplicity:

{\begin{assumption}
The system~$\Sigma_{\theta}$ has the same number of inputs and outputs, that is, $n_y = n_w = m$ for a positive integer~$m$. 
\end{assumption}}

This assumption simplifies the notation and is not restrictive because we can insert the input and output matrices with zero columns and rows to realize $n_w = n_y$, without affecting the robust stability notions we shall discuss below (see also, e.g.,~\cite{Colombino2016a}). We then consider the situation in which the open-loop system~$\Sigma_{\theta}$ is closed with the relationship
\begin{equation}\label{eq:interconnection}
w = \Delta y, 
\end{equation}
where $\Delta \in \bm \Delta \subset [0, \infty)^{m\times m}$ represents a static uncertainty matrix. In this section, we are interested in the stability of the closed-loop system arising from the interconnection, that is, the internal stability of the system
\begin{equation}\label{eq:closedLoop}
\frac{dx}{dt} = \left(A(\theta) + B(\theta)\Delta C(\theta)\right)x. 
\end{equation}
To quantify the robust stability of this closed-loop system, let us introduce the quantity 
\begin{equation*}
\eta(\theta) = \sup_{\Delta\in \bm \Delta,\,\norm{\Delta}\leq {\epsilon}} \lambda_{\max}(A(\theta)+B(\theta) \Delta C(\theta)), 
\end{equation*}
{where $\epsilon > 0$ represents the maximum size of the uncertainty matrix~$\Delta$.} In this context, we consider the following robust stabilization problem:
\begin{subequations}\label{eq:strucutredStabilization}
\begin{align}
\minimize_{\theta\in\Theta}\ \ \ & L(\theta)
\\
\subjectto\ \ \,\,
& \eta(\theta) < -\gamma,\label{eq:strucutredStabilization:eta}
\end{align}
\end{subequations}
where $\gamma > 0$ denotes the desired exponential decay rate for the closed-loop system~\eqref{eq:closedLoop}.

Following the formulation in~\cite{Colombino2016a}, this paper focuses on the structural uncertainties belonging to  
\begin{equation*}
\begin{aligned}
\mathcal D &= \{
\diag (\Delta_1, \dotsc, \Delta_\phi, \delta_{\phi+1}, \dotsc, \delta_{\phi+\sigma})\mid 
\\
&\quad \quad \quad  \Delta_k \in \mathbb{R}_+^{m_k\times m_k},\ k=1, \dotsc, \phi, 
\\
&\quad \quad \quad \delta_k \geq  0 ,\ k = \phi+1, \dotsc, \phi+\sigma
\}\subset \mathbb{R}^{m\times m}. 
\end{aligned} 
\end{equation*}
Then, the following theorem shows that we can solve the robust stabilization problem~\eqref{eq:strucutredStabilization} by geometric programming. 

\begin{theorem}\label{thm:mu}
Define the set
\begin{equation*}
\begin{aligned}
\mathcal P &= \{
\diag(\pi_1 I_{m_1}, \dotsc, \pi_\phi I_{m_\phi}, \pi_{\phi+1}, \dotsc, {\pi}_{\phi+\sigma})\mid \\
&\quad\quad \quad{\pi}_k >0,\ k=1, \dotsc, \phi + \sigma
\}\subset 
\mathbb{R}^{m\times m}. 
\end{aligned}
\end{equation*}
Then, the solution of the robust stabilization problem~\eqref{eq:strucutredStabilization} is given by the following geometric program: 
\begin{subequations}\label{eq:muopt}
\begin{align}
\minimize_{\mathclap{\substack{\theta\in\mathbb R^\nsubtheta_{++},\, \Pi \in \mathcal P ,\\u, v\in\mathbb{R}^m_{++},\\\xi, \zeta\in\mathbb{R}^{n_x}_{++}}}}\ \ & \tilde L(\theta)
\\
\subjectto\ \,\,& {\sqrt \epsilon}D_v^{-1}\Pi^{1/2}  C(\theta) \xi  <  \onev,\label{eq:muoptconstfirst}
\\
&D_\xi^{-1} R(\theta)^{-1}(\tilde A(\theta)\xi + \gamma \xi +{\sqrt \epsilon} B(\theta)\Pi^{-1/2} u) < \onev, \!\!\!
\\
& {\sqrt \epsilon} D_u^{-1}\Pi^{-1/2} B(\theta)^\top \zeta < \onev,
\\
& D_\zeta^{-1}R(\theta)^{-1}(\tilde A(\theta)^\top \zeta +\gamma \zeta +  {\sqrt \epsilon}  C(\theta)^\top \Pi^{1/2} v) < \onev,\!\!\!
\\
&f_i(\theta) \leq 1, \quad i=1, \dotsc, p.  \label{eq:muoptconstlast}
\end{align}
\end{subequations}
\afterequation
\end{theorem}

In order to prove Theorem~\ref{thm:mu}, we present the following proposition.

\begin{proposition}\label{prop:robustStability}
Consider the positive linear system~$\Sigma$ given by~\eqref{eq:lti}. Let $\gamma > 0$. The following two conditions are equivalent: 
\begin{enumerate}
\item The following inequality holds true: 
\begin{equation}\label{eq:leq-gamma}
\sup_{\Delta\in \bm \Delta,\,\norm{\Delta}\leq {\epsilon}} \lambda_{\max}(F+G \Delta H)< -\gamma. 
\end{equation}
\item There exist positive vectors~$u, v \in \mathbb{R}^m$ and~$\xi, \zeta\in\mathbb{R}^{n_x}$ as well as a matrix~$\Pi \in {\mathcal P}$ such that the following inequalities hold true:
\begin{equation}
\begin{aligned}
&{\sqrt \epsilon}\Pi^{1/2} H \xi  < v, \label{eq:gather1}
\\
&(F+\gamma I)\xi +{\sqrt \epsilon} G\Pi^{-1/2} u < 0, 
\\
&{\sqrt \epsilon}\Pi^{-1/2} G^\top \zeta < u,
\\
&(F^\top + \gamma I) \zeta + {\sqrt \epsilon}H^\top \Pi^{1/2} v < 0.
\end{aligned}
\end{equation}
\afterequation
\end{enumerate}
\end{proposition}

\begin{IEEEproof}
Let us prove the necessity. Assume that inequality~\eqref{eq:leq-gamma} holds true. Then, the system
\begin{equation*}
\Sigma_\gamma\colon
\begin{cases}
\,\dfrac{dx}{dt} = (F+\gamma I)x + {\sqrt \epsilon} Gw, 
\\
\,y = {\sqrt \epsilon} Hx,\end{cases}  
\end{equation*}
with the feedback~\eqref{eq:interconnection} is internally stable for all $\Delta\in \bm\Delta$ satisfying $\norm{\Delta}\leq 1$. Let $\hat M_\gamma(s)$ denote the transfer function of the system~$\Sigma_\gamma$. Then, by~\cite[Theorem~10]{Colombino2016a}, there exists $\Pi \in \mathcal P$ such that $\norm{\Pi^{1/2} \hat M_\gamma(0)\Pi^{-1/2}} < 1$. Therefore, Lemma~\ref{lem:normmat} shows the existence of positive vectors~$u, v\in \mathbb{R}^{m}$ such that
\begin{equation*}
\begin{aligned}
&{-\Pi^{1/2}} {\sqrt \epsilon} H(F+\gamma I)^{-1}{\sqrt \epsilon} G\Pi^{-1/2}u < v, 
\\
&{- \Pi^{-1/2}} {\sqrt \epsilon} G^\top (F^\top+\gamma I)^{-1} {\sqrt \epsilon} H^\top  \Pi^{1/2}v  < u. 
\end{aligned}
\end{equation*}
In the same way as in the proof of Proposition~\ref{prop:Hinf}, applying Lemma~\ref{lem:key} to these inequalities shows the existence of positive vectors~$\xi, \zeta \in \mathbb{R}^{n_x}$ satisfying the inequalities in~\eqref{eq:gather1}, as desired. The proof of sufficiency is omitted.
\end{IEEEproof}

Let us prove Theorem~\ref{thm:mu}. 

\begin{proofof}{Theorem~\ref{thm:mu}}
Proposition~\ref{prop:robustStability} implies that the solution of the robust stabilization problem~\eqref{eq:strucutredStabilization} is given by the solution of the following optimization problem: 
\begin{equation*}
\begin{aligned}
\minimize_{\mathclap{\substack{\theta\in\Theta,\,\Pi \in \mathcal P ,\\u, v \in \mathbb{R}^m_{++},\,\xi, \zeta\in\mathbb{R}^{n_x}_{++}}}}\ \ \ \  & L(\theta)
\\
\subjectto\ \ \ \,\, &  {\sqrt \epsilon}\Pi^{1/2} C(\theta) \xi  <  v,
\\
& (A(\theta) + \gamma I)\xi + {\sqrt \epsilon} B(\theta)\Pi^{-1/2} u < 0,
\\
& {\sqrt \epsilon}\Pi^{-1/2}  B(\theta)^\top \zeta < u,
\\
& (A(\theta)+ \gamma I)^\top \zeta + {\sqrt \epsilon} C(\theta)^\top \Pi^{1/2} v < 0. 
\end{aligned}
\end{equation*}
A simple algebraic manipulation reduces this optimization problem to the optimization problem~\eqref{eq:muopt}, which is indeed a geometric program by Assumptions~\ref{asm:posy} and~\ref{asm:costPosy} as well as the fact that $\Pi$ is a diagonal matrix whose diagonals are monomials with respect to the variables $\pi_k$. The further details of the proof are omitted.
\end{proofof}

Finally, as a direct corollary of Theorem~\ref{thm:mu}, we below present a geometric program for identifying the maximum allowable size of the uncertainty matrix~$\Delta$ for the robust stabilization problem~\eqref{eq:strucutredStabilization} to be feasible. 

\begin{corollary}\label{cor:}
The robust stabilization problem~\eqref{eq:strucutredStabilization} is solvable for all $\epsilon \in [0, \epsilon^\star]$, where $\epsilon = \epsilon^\star$ solves the following geometric program: 
\begin{equation*}
\begin{aligned}
\minimize_{\mathclap{\substack{\theta\in\mathbb R^\nsubtheta_{++},\, \Pi \in \mathcal P ,\\u, v\in\mathbb{R}^m_{++},\,\xi, \zeta\in\mathbb{R}^{n_x}_{++},\,\epsilon> 0}}}\ \ \ \ \ \  \  & 1/\epsilon
\\
\subjectto\ \ \ \ \ \ \,\,& \mbox{\eqref{eq:muoptconstfirst}--\eqref{eq:muoptconstlast}}.
\end{aligned}
\end{equation*}
\afterequation
\end{corollary}

\section{Time-delay systems} \label{sec:delay}

In the previous sections, we have presented geometric programming-based frameworks for efficiently solving various classes of norm-constrained {parameter }optimization problems for positive linear systems. The aim of this section is to extend the frameworks to delayed positive linear systems~\cite{Haddad2004}. Let us consider the following parametrized positive linear system with time-delays: 
\begin{equation*}
\Sigma_{d, \theta}:
\begin{cases}
\,\dfrac{dx}{dt} = A(\theta)x(t) +A_d(\theta)x(t-h) + B(\theta)w(t), 
\\
\,y = C(\theta)x(t) + C_d(\theta)x(t-h), 
\\
\,x\arrowvert_{[-h, 0]} = \phi\in \mathcal C([-h, 0], \mathbb{R}_+^{n_x}), 
\end{cases}  
\end{equation*}
where $h>0$ represents a constant delay and  $\mathcal C([-h, 0], \mathbb{R}_+^{n_x})$ denotes the set of~$\mathbb{R}_+^{n_x}$-valued continuous functions defined on the interval $[-h, 0]$. We denote the solutions of the system~$\Sigma_{d, \theta}$ with the initial condition~$\phi$ and the disturbance signal~$w$ by~$x(t; \phi, w)$ and~$y(t; \phi, w)$, when we need to emphasize their dependence on $\phi$ and~$w$.  We suppose that, for all $\theta \in \Theta$, the matrix~$A(\theta)$ is Metzler and the matrices~$A_d(\theta)$, $B(\theta)$, $C(\theta)$, and~$C_d(\theta)$ are nonnegative. This guarantees~\cite{Liu2011d} that the system~$\Sigma_{d, \theta}$ is internally positive, {that is}, the values of the state~$x(t)$ and output $y(t)$ remain nonnegative at every time instant~$t$ if $\phi(t)\geq 0$ for all $t \in [-h, 0]$ and~$w(t)\geq 0$ for all $t\geq 0$. 

We are concerned with the following three quantities on the delayed positive linear system~$\Sigma_{d, \theta}$. The first one is the exponential decay rate 
defined by 
\begin{equation*}
\rho_\theta = -\sup_{\phi\in \mathcal C([-h, 0], \mathbb{R}_+^n)}\limsup_{t\to\infty}\frac{\log \norm{x(t; \phi, 0)}}{t}. 
\end{equation*}
The second quantity of interest is the $\mathcal L^1$-gain of the system~\cite{Briat2012c}, {\cite{Desor1975}}, \cite{Zhu2017}. Assume that $\rho_\theta > 0$. For a positive integer~$n$, let $\mathcal L^1(\mathbb{R}^{n}_+) = \{ f\colon [0, \infty) \to \mathbb{R}^{n} \mid \int_0^\infty \norm{f(t)}_1\,dt < \infty \}$ denote the space of Lebesgue-measurable integrable functions equipped with the norm $\norm{f}_1 = \int_0^\infty \norm{f(t)}_1\,dt$, where $\norm{f(t)}_1$ denotes the $1$-norm of the vector~$f(t)$. Then, we define the $\mathcal L^1$-gain of~$\Sigma_{d, \theta}$ by 
\begin{equation*}
\norm{\Sigma_{d, \theta}}_{\mathcal L^1} = \sup_{w\in \mathcal L^1(\mathbb{R}^{n_w}_{+}) \backslash \{0\}} \frac{\norm{y(\cdot; 0, w)}_{ 1}}{\norm{w}_{1}}.
\end{equation*}
The third and last quantity of our interest is the $\mathcal L^\infty$-gain~{\cite{Desor1975}}, \cite{Shen2015}. Let $\mathcal L^\infty(\mathbb{R}^{n}_+) = \{ f\colon [0, \infty) \to \mathbb{R}^{n} \mid \esssup_{t\geq 0} \norm{f(t)}_\infty  < \infty \}$ denote the space of~$\mathbb{R}^{n}_+$-valued essentially bounded Lebesgue-measurable functions equipped with the norm $\norm{f}_\infty = \esssup_{t\geq 0} \norm{f(t)}_\infty$, where $\norm{f(t)}_\infty$ denotes the $\infty$-norm. Then, we define the $\mathcal L^\infty$-gain of~$\Sigma_{d, \theta}$ by
\begin{equation*}
\norm{\Sigma_{d, \theta}}_{\mathcal L^\infty} = \sup_{w\in \mathcal L^\infty(\mathbb{R}^{n_w}_{+}) \backslash \{0\}} \frac{\norm{y(\cdot; 0, w)}_{\infty}}{\norm{w}_{\infty}}.
\end{equation*}
Then, the {parameter }optimization problem that we study in this section is stated as follows: 
\begin{subequations}\label{eq:delayMixedOptimization}
\begin{align}
\minimize_{\theta\in\Theta}\ \ \ & L(\theta)
\\
\subjectto\ \ \,\,& 
\rho_\theta > 0, \label{eq:requirement:rho}
\\
&\beta(\rho_\theta, \norm{\Sigma_\theta}_{\mathcal L^1}, \norm{\Sigma_\theta}_{\mathcal L^\infty}) < \gamma,\label{eq:requirement:rhoL1Linf}
\end{align}
\end{subequations}
where $\beta\colon \mathbb{R}_{++}^3\to\mathbb{R}_{++}$ is a function representing the trade-off between the exponential decay rate, $\mathcal L^1$-gain, and~$\mathcal L^\infty$-gain of the system. 

Let us place the following assumption, which corresponds to Assumptions~\ref{asm:posy} and~\ref{asm:alpha} in the delay-free case. 

\begin{assumption}\label{asm:posy:d}
The following conditions hold true: 
\begin{enumerate}
\item 
There exists a diagonal matrix function $R(\theta)$
having monomial diagonals such that each entry of the matrix
\begin{equation*}
\tilde A_d(\theta) = A(\theta) + A_d(\theta) + R(\theta)
\end{equation*}
is either  a posynomial of~$\theta$ or zero. 
\item  
Each entry of the matrices $B(\theta)$, $B_d(\theta)$, $C(\theta)$, and~$C_d(\theta)$ is either a posynomial of~$\theta$ or zero. 
\item The function $\beta$ is a posynomial, nonincreasing with respect  to the first variable, and nondecreasing with respect to the left two variables. \label{item:asm:beta}
\end{enumerate}
\afterequation
\end{assumption}

Under these assumptions, the following theorem shows that the mixed-constraint optimization problem~\eqref{eq:delayMixedOptimization} can be solved by convex optimization. 

\begin{theorem}\label{thm:mixed:delay}
Let $\gamma >0$ be given. Define the function~$g$ by 
\begin{equation*}
g(\rho) = e^{\rho h} - 1
\end{equation*}
for $\rho > 0$. The solution of the {parameter }optimization problem~\eqref{eq:delayMixedOptimization} is given by the solution of the following optimization problem: 
\begin{subequations}\label{eq:optL1inf:pre}
\begin{align}
\minimize_{\mathclap{\substack{\theta\in \mathbb R^\nsubtheta_{++},\\\xi,\,u,\, v\in\mathbb{R}^{n_x}_{++},\\\rho,\,\gamma_1,\,\gamma_\infty > 0}}}\ \ \ & \tilde L(\theta)
\\
\subjectto\ \ \,\,&
\gamma^{-1}\beta(\rho, \gamma_1, \gamma_\infty) < 1, \label{eq:beta<gamma}
\\
&D_\xi^{-1}R(\theta)^{-1}\bigl(\tilde A_d(\theta)+\rho I + g(\rho)A_d(\theta)\bigr)\xi < \onev, \label{eq:delay:stability}
\\
&D_u^{-1}R(\theta)^{-1}\bigl(\tilde A_d(\theta)^\top u+(C(\theta)+C_d(\theta))^\top \onev\bigr) <  \onev, \label{eq:delay:L1:1}
\\
& \gamma_1^{-1}B(\theta)^\top u < \onev,\label{eq:delay:L1:2}
\\
& D_v^{-1} R(\theta)^{-1} (\tilde A_d(\theta)v+B(\theta)\onev)< \onev,\label{eq:delay:Linf:1} 
\\
& \gamma_\infty^{-1}(C(\theta)+C_d(\theta))v < \onev, \label{eq:delay:Linf:2}
\\
&f_i(\theta) \leq 1, \quad i=1, \dotsc, p. \label{eq:delay:f_i}
\end{align}
\end{subequations}
Moreover, this optimization problem reduces to a convex optimization problem by logarithmic variable transformations of the form~\eqref{eq:transformation}.
\end{theorem}

\begin{remark}
The optimization problem~\eqref{eq:optL1inf:pre} is not a geometric program because the function $g$ appearing in the constraint~\eqref{eq:delay:stability} is not a posynomial. However, as stated in Theorem~\ref{thm:mixed:delay} shall be shown below in the proof of the theorem, the optimization problem can be still reduced to a convex optimization problem. 
\end{remark}

\begin{IEEEproof}[Proof of Theorem~\ref{thm:mixed:delay}]
Assumptions~\ref{asm:costPosy} and~\ref{asm:posy:d} show that the optimization problem~\eqref{eq:optL1inf:pre} is a geometric program if the function~$g$ was a posynomial. However, because $g$ is not a posynomial, the optimization problem~\eqref{eq:optL1inf:pre} is not a geometric program. However, the function~$g$ is a limit of the sequence of posynomials~$\{g_k\}_{k=1}^\infty$ given by $g_k(\rho) = \sum_{\ell=1}^{k}({\rho h})^\ell/\ell!$. Therefore, logarithmic variable transformations of the form~\eqref{eq:transformation} in fact convert the optimization problem~\eqref{eq:optL1inf:pre} into a convex optimization problem (see~\cite[Section~7.1]{Boyd2007} for further details).

As in the proof of Corollary~\ref{cor:mixedH2Hinfty}, we need to show that $\theta \in \Theta$ satisfies inequalities~\eqref{eq:requirement:rho} and~\eqref{eq:requirement:rhoL1Linf} if and only if there exist positive vectors~$\xi, u, v\in \mathbb{R}^{n_x}$ and positive numbers~$\rho, \gamma_1, \gamma_\infty$ satisfying  constraints~\eqref{eq:beta<gamma}--\eqref{eq:delay:Linf:2}. 

In this proof, we only show the sufficiency. Suppose the existence of positive vectors~$\xi, u, v\in \mathbb{R}^{n_x}$ and positive numbers~$\rho, \gamma_1, \gamma_\infty$ satisfying \eqref{eq:beta<gamma}--\eqref{eq:delay:Linf:2}. By the monotonicity property of the function~$\beta$ (see Assumption~\ref{asm:posy:d}.\ref{item:asm:beta}) and inequality~\eqref{eq:beta<gamma}, it is sufficient to show the following inequalities 
\begin{align}
\rho_\theta &> \rho,\label{eq:rhotheta<rho} 
\\
\norm{\Sigma_{d, \theta}}_{\mathcal{L}^1} &< \gamma_1, \label{eq:norm<gamma1} 
\\
\norm{\Sigma_{d, \theta}}_{\mathcal{L}^\infty} &< \gamma_\infty.  \label{eq:norm<gammainfinity} 
\end{align}
Let us first show \eqref{eq:rhotheta<rho}. Let $\phi \in \mathcal C([-h, 0], \mathbb{R}_+^{n_x})$ be arbitrary. Since inequality~\eqref{eq:delay:stability} implies $(A(\theta)+\rho I + e^{{\rho h}}A_d(\theta))\xi < 0$, Lemma~\ref{lem:PF} shows that the matrix~$A(\theta)+\rho I + e^{{\rho h}}A_d(\theta)$ is Hurwitz stable. Therefore, Theorem~3.1 in~\cite{Ngoc2013} shows that the solution~$y$ of the following delayed positive linear system
\begin{equation*}
\frac{d\tilde x}{dt} = (A(\theta) + \rho I)\tilde x + e^{{\rho h}}A_d(\theta)\tilde x(t-h)
\end{equation*}
converges to zero exponentially fast. On the other hand, the function $\tilde x(t) = e^{\rho t}x(t;\phi, 0)$ satisfies this differential equation for $t\geq h$. Therefore, we conclude that the function $x(\cdot; \phi, 0)$ converges to zero exponentially fast with its rate being greater than $\rho$, as desired. 
We then show inequalities~\eqref{eq:norm<gamma1} and~\eqref{eq:norm<gammainfinity}. Inequalities~\eqref{eq:delay:L1:1} and \eqref{eq:delay:L1:2} show $(A(\theta)+A_d(\theta))^\top u + (C(\theta)+C_d(\theta))^\top \onev < 0$ and~$B^\top(\theta) u -\gamma_1 \onev < 0$. These inequalities and Lemma~2 in~\cite{Zhu2017} show \eqref{eq:norm<gamma1}. In a similar manner, Theorem~2 in~\cite{Shen2015} shows that inequalities \eqref{eq:delay:Linf:1} and \eqref{eq:delay:Linf:2} imply \eqref{eq:norm<gammainfinity}. This completes the proof of the theorem.  
\end{IEEEproof}

\section{Example: dynamical buffer networks} \label{sec:simulations2}

In this section, we illustrate the theoretical results presented in the previous sections. Let $\mathcal G$ be a weighted and directed graph with the node set~$\mathcal V = \{1, \dotsc, N\}$ and edge set~$\mathcal E = \{e_1, \dotsc, e_M\} \subset \{1, \dotsc, N\} \times \{1, \dotsc, N\}$, respectively. For each edge $e_\ell$ we use the notation $e_\ell = (e_\ell(1), e_\ell(2))$, where the nodes $e_\ell(1)$ and $e_\ell(2)$ denote the origin and the destination of the edge, respectively. Since the graph~$\mathcal G$ is weighted, a positive and fixed weight $w_{e_\ell}$ is assigned on an edge~$e_\ell$. By abusing the notation, we often write the weight $w_{e_\ell}$ as $w_{e_\ell(1) e_\ell(2)}$. Therefore, the weight of an edge $(i, j)$ is denoted by $w_{ij}$. We define the adjacency matrix~$A_{\mathcal G}\in \mathbb{R}^{N\times N}$ of the graph~$\mathcal G$ by
\begin{equation*}
[A_{\mathcal G}]_{ij} = \begin{cases}
w_{ji},&\mbox{if $(j, i)\in \mathcal E$, }
\\
0,&\mbox{othereise.}
\end{cases}
\end{equation*}
Also, let us define the set of in-neighborhood of node~$i$ as $\mathcal N^{\textrm{in}}_i = \{ j \in \mathcal V: (j, i) \in \mathcal E\}$. Similarly, we define the set of out-neighborhood of node~$i$ as $\mathcal N^{\textrm{out}}_i = \{ j \in \mathcal V: (i, j) \in \mathcal E\}$.

We assume that the network~$\mathcal G$ contains at least one origin (i.e., a node having an empty in-neighborhood) and at least one destination (i.e., a node having an empty out-neighborhood). Let $\mathcal V_o$ and~$\mathcal V_d$ denote the set of origins and destinations of the network, respectively. Without loss of generality, we assume that $\mathcal V_o = \{1, \dotsc, \abs{\mathcal V_o}\}$, where $\abs{\mathcal V_o}$ denotes the size of the set~${\mathcal V_o}$. We allow the network to have multiple origins and/or neighbors. Then, we consider a dynamical buffer network described by the following set of differential equations (see, e.g.,~\cite{Rantzer2018}):
\begin{equation}\label{eq:dxidt/flow}
\frac{dx_i}{dt} = 
\begin{cases}
\displaystyle 
{f^{\textrm{in}}_i - \sum_{j \in \mathcal N^{\textrm{out}}_i} u_{ij}}, 
& \mbox{if $i \in \mathcal V_o$, }
\\\displaystyle
 \sum_{j \in \mathcal N^{\textrm{in}}_i} u_{ji}
- \sum_{j \in \mathcal N^{\textrm{out}}_i} u_{ij},  
&\displaystyle \mbox{if $i \notin \mathcal V_o \cup \mathcal V_d$,}
\\\displaystyle
{\sum_{j \in \mathcal N^{\textrm{in}}_i} u_{ji} - f^\textrm{out}_i}, 
& \mbox{if $i \in \mathcal V_d$,}
\end{cases}
\end{equation}
where $x_i$ represents the buffer content of node~$i$, $u_{ij}$ represents the volume of flow from node $i$ to $j$, $f^\textrm{in}_i$ ($i\in\mathcal V_o$) describes the effect of local production or an external disturbance, and $f^\textrm{out}_i$ ($i\in\mathcal V_d$) describes the decay of the buffer content at destination nodes. The flows are assumed to be in the following linear form: 
\begin{equation*}
f^{\textrm{out}}_i = \phi_i x_i,\ u_{ij} = \psi_i w_{ij} x_i, 
\end{equation*}
where $\phi_i > 0$ ($i \in \mathcal V_d$) and ${\psi_i>0}$ ($i \in \mathcal V \backslash \mathcal V_d$) are constants dependent on node~$i$. For convenience of notation, we set $\phi_i = 0$ for all $i \in \mathcal V \backslash \mathcal V_d$ and $\psi_i = 0$ for all $i \in \mathcal V_d$. Also, let us set the measurement output of the system as
\begin{equation}\label{eq:y=:buffer}
y = \begin{bmatrix}
x \\ \alpha u
\end{bmatrix}, 
\end{equation}
where $\alpha > 0$ is a weight constant and the $M$-dimensional vector~$u$ is obtained by vertically stacking the flows $u_{ij}$ as $u = [ u_{e_1(1)e_1(2)} \  \cdots \  u_{e_M(1)e_M(2)} ]^\top$. Let us denote the dynamical system~\eqref{eq:dxidt/flow} and~\eqref{eq:y=:buffer} by $\Sigma_{\phi, \psi}$, which we can rewrite as 
\begin{equation*}
\Sigma_{\phi, \psi}: \begin{cases}
\displaystyle \frac{dx}{dt} = \left(
A_{\mathcal G}\Psi - \diag(\onev^\top A_{\mathcal G}\Psi) - \Phi
\right)x + \begin{bmatrix}
I_{\abs{\mathcal V_o}} \\ O_{n-\abs{\mathcal V_o}, n}
\end{bmatrix}f^\textrm{in},
\\
\displaystyle y = \begin{bmatrix}
I_n \\ \alpha H \Psi
\end{bmatrix} x, 
\end{cases}
\end{equation*}
where $f^\textrm{in} = [f^\textrm{in}_1\ \cdots \ f^\textrm{in}_{\abs{\mathcal V_o}}]^\top$, $\Psi = \diag(\psi_1, \dotsc,\psi_N)$, $\Phi = \diag(\phi_1, \dotsc, \phi_N)$, and the matrix~$H \in \mathbb{R}^{M\times N}$ is defined by 
\begin{equation*}
H_{\ell i} = \begin{cases}
w_{e_\ell}, & \mbox{if $i = e_\ell(1)$,}
\\
0,&\mbox{otherwise,}
\end{cases}
\end{equation*}
for all $\ell \in \{1, \dotsc, M \}$ and~$i \in \{ 1, \dotsc, N \}$. 

In this example, we study the problem of tuning the local constants~$\phi_i$ and~$\psi_i$ for achieving a small $H^\infty$ norm of the dynamical buffer network $\Sigma_{\phi, \psi}$. Let us introduce the variables~$\phi = (\phi_i)_{i \in \mathcal V_d}$ and $\psi = (\psi_i)_{i\in\mathcal V\backslash \mathcal V_d}$. We measure the cost for tuning the system by the sum
\begin{equation}\label{eq:L(theta)}
L(\phi, \psi) = \sum_{i \in \mathcal V_d} \phi_i + \sum_{i \in \mathcal V \backslash \mathcal V_d}\psi_i.  
\end{equation}
We further allow the following forms of upper-bounds on the parameters to be tuned: 
\begin{equation}\label{eq:box:phipsi}
\phi_i \leq \bar \phi_i,\ \psi_i \leq \bar \psi_i
\end{equation}
for positive constants $\bar \phi_i$ and $\psi_i$, which may arise from physical restrictions. We can now formulate our optimization problem as follows.

\begin{problem}[Buffer network optimization]
Let $\gamma > 0$ be given. Find the set of parameters~$\phi$ and~$\psi$ satisfying the constraints in~\eqref{eq:box:phipsi} as well as the $H^\infty$ norm-constraint~$\norm{\Sigma_{\phi,\psi}}_\infty < \gamma$, while the cost~$L(\phi, \psi)$ is minimized.
\end{problem}

Let us show that the buffer network optimization problem can be solved by geometric programming. It is easy to see that the system~$\Sigma_{\phi, \psi}$ satisfies Assumption~\ref{asm:posy}.\ref{asm:posy:BC}). In order to show that Assumption~\ref{asm:posy}.\ref{asm:posy:A}) is satisfied, we define the matrix~$R(\theta) = \diag(\onev^\top A_{\mathcal G}\Psi) + \Phi$. Then, each entry of the matrix~$\tilde A(\theta) = A(\theta) + R(\theta) = A_{\mathcal G}\Psi$ is either a posynomial in the variables $\phi$ and~$\psi$ or zero. Moreover, $R(\theta)$ is a diagonal matrix and has the monomial diagonals:
\begin{equation*}
[R(\theta)]_{ii} = 
\begin{cases}
\psi_i \sum_{j \in \mathcal N^{\textrm{out}}_i} w_{ij}, & \mbox{if $i \in \mathcal V \backslash  \mathcal V_d$,}
\\
\phi_i, & \mbox{otherwise}.
\end{cases}
\end{equation*}
Therefore, Assumption~\ref{asm:posy}.\ref{asm:posy:A}) is satisfied as well.  Also, it is trivial to see that the cost function~\eqref{eq:L(theta)} is a posynomial in the variables and, therefore, satisfies Assumption~\ref{asm:costPosy}.\ref{ams:item:L}). Finally,  Assumption~\ref{asm:costPosy}.\ref{ams:item:Theta}) holds true because one can rewrite the constraints~\eqref{eq:box:phipsi} in terms of posynomials as $\bar \phi_i^{-1}\phi_i \leq 1$ and $\bar \psi_i^{-1}\psi_i \leq 1$. Therefore, we can use Theorem~\ref{thm:Hinf} to solve the buffer network optimization problem via geometric programming.

\begin{figure}[tb]
\centering
\includegraphics[width=.65\linewidth]{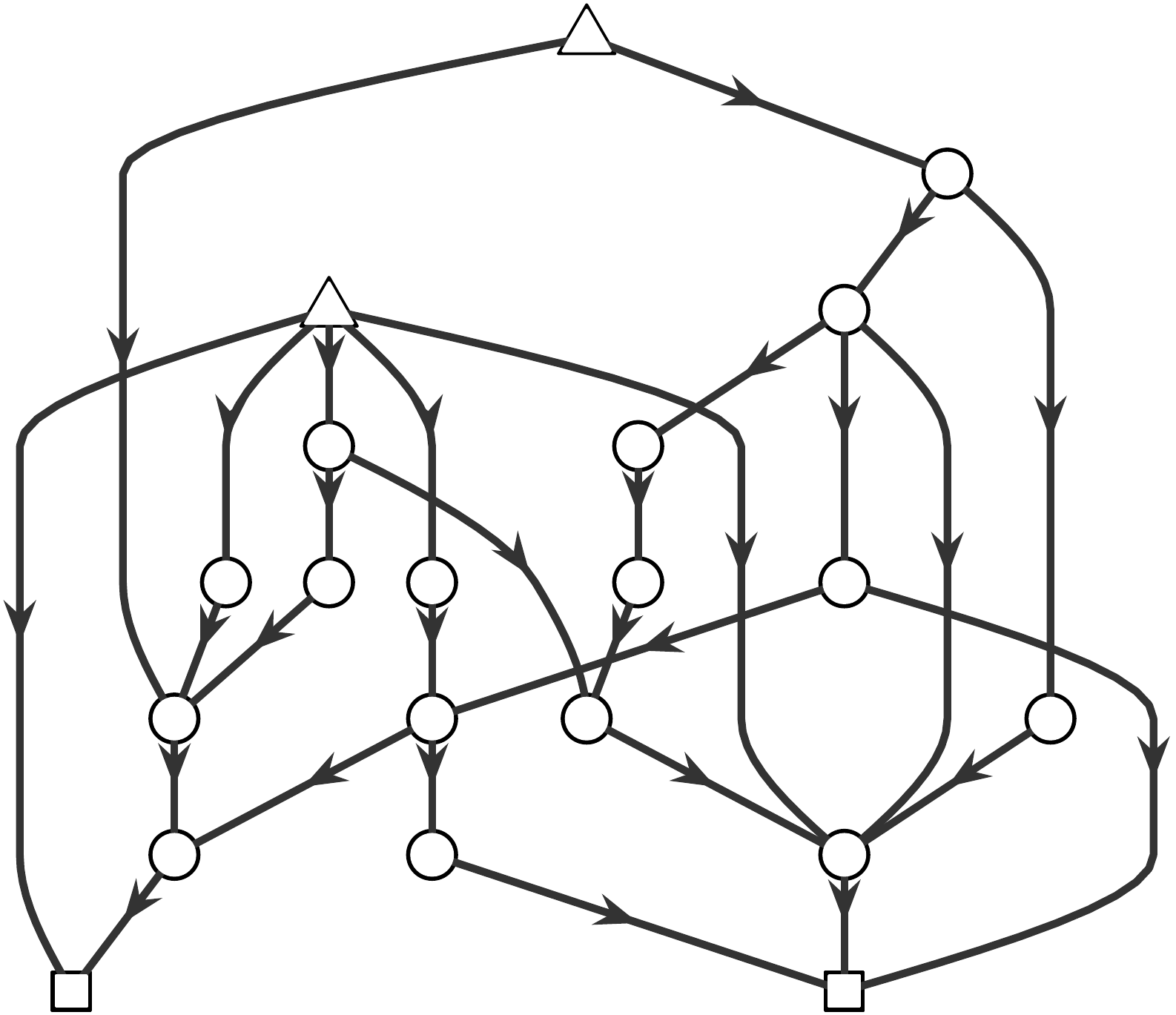}
\caption{A directed acyclic graph. Triangles: origins. Squares: destinations.}
\label{fig:flow}
\end{figure}

For numerical simulations, let us consider a synthetic directed acyclic graph shown in Fig.~\ref{fig:flow}. The graph has two origins (indicated by triangles) and two destinations (indicated by squares). We assume that the weights of edges originating from a node are equal and sums to one. Therefore, we set $w_{ij} = 1/\abs{\mathcal N_i^{\textrm{out}}}$ for all node~$i$. Also, let us set $\bar \phi_i = \bar \psi_i = 5$ for all nodes and use the weight $\alpha = 1/10$ in the measurement output~\eqref{eq:y=:buffer}. Using an $H^\infty$ norm-counterpart of Remark~\ref{rmk:}, we first identify the minimum achievable $H^\infty$ norm of the system as $\gamma^\star = 0.388$. Then, for various values of~$\gamma$ within the interval $[\gamma^\star, 4\gamma^\star]$, we solve the buffer network optimization problem and obtain the optimal values of the local parameters~$\phi$ and~$\psi$. We show the values of the optimal cost~$L$ for various values of~$\gamma$. For the cases when $\gamma = 1.5\gamma^\star$, $2\gamma^\star$, and~$4\gamma^\star$, we illustrate the obtained values of the constants~$\phi$ and~$\psi$ in Fig.~\ref{fig:flow:investments}.

\begin{figure}[tb]
\centering
\includegraphics[width=.98\linewidth]{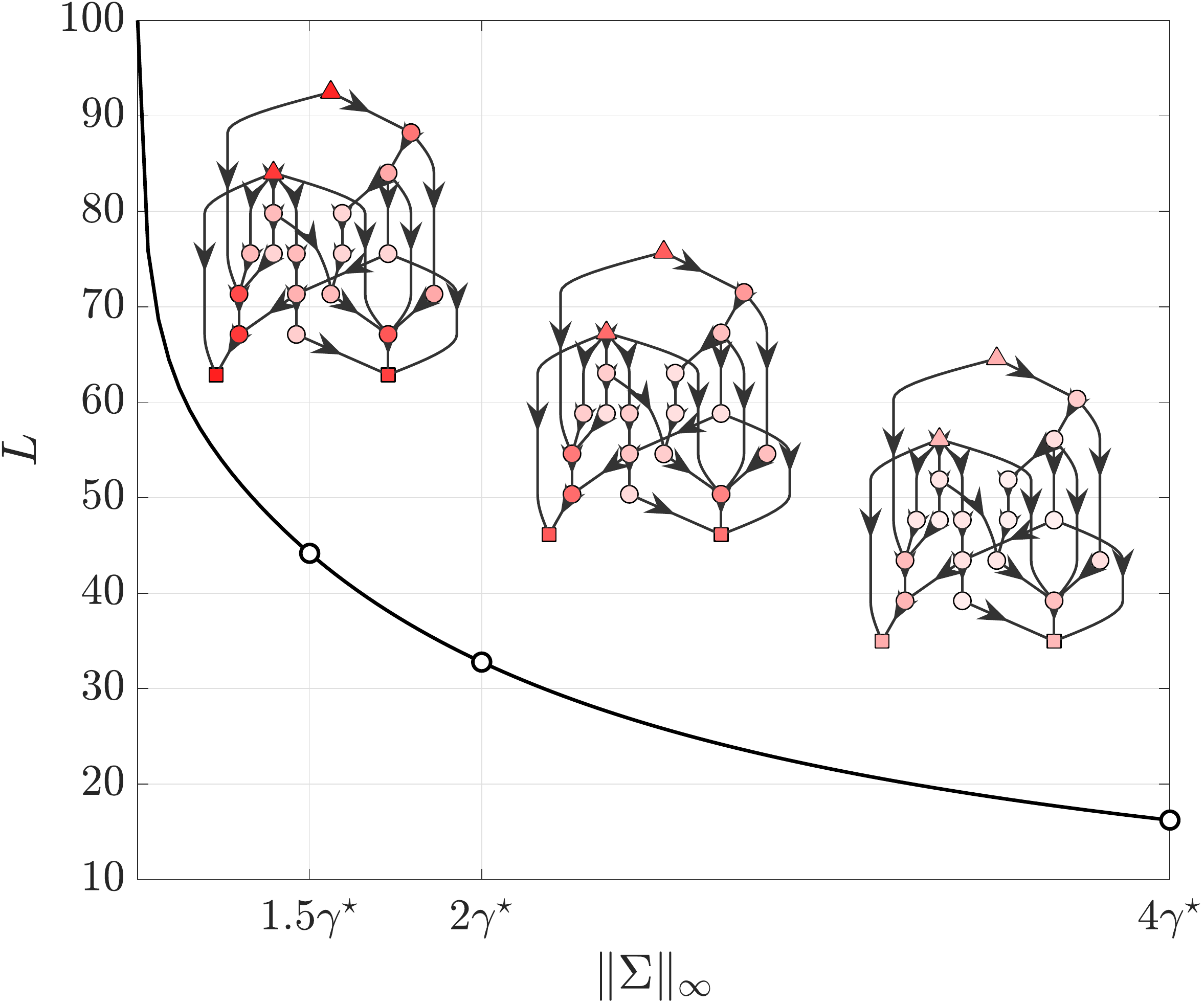}
\caption{The amount of optimal investments for various values of~$\gamma$. The darker the colors, the bigger the optimized parameters $\phi_i$ and~$\psi_i$.}
\label{fig:flow:investments}
\end{figure}

\section{Example: networked epidemics} \label{sec:simulations}

In this section, we consider the Susceptible-Infected-Susceptible (SIS) model for describing networked epidemic processes taking place in human and animal social networks~\cite{Nowzari2015a,Pastor-Satorras2015a}. In the SIS model, at a given (continuous) time~$t \geq 0$, each node can be in one of two possible states: {\it susceptible} or {\it infected}. When a node~$i$ is infected, it can randomly {transit} to the susceptible state with an instantaneous rate~$\delta_i > 0$, called the {\it recovery rate} of node~$i$. On the other hand, an infected node~$j$ can infect node~$i$ with the instantaneous rate~$\beta_iw_{ij}$, where $\beta_i > 0$ is called the {\it infection rate} of node $i$. The SIS model is a Markov process having a total of~$2^N$ possible states~\cite{VanMieghem2009a} (two states per node). 

Throughout this section, we consider the situation where the connectivity of the network is not completely known, as is often the case in practice~\cite{Holme2017}. In this paper, let us model this uncertainty as an additive uncertainty in the weights of edges, {that is}, let us assume that the adjacency matrix of the graph takes the form
\begin{equation*}
A_{\mathcal G} = A_{\bar {\mathcal G}} + A_{\Delta \mathcal G}, 
\end{equation*}
where $A_{\bar {\mathcal G}}$ denotes the adjacency matrix of the nominal (weighted) network~$\bar {\mathcal G}$ and~$A_{\Delta\mathcal G}$ denotes a nonnegative matrix representing the uncertainty. For simplicity, we assume that only a bound on the norm of the uncertainty $A_{\Delta\mathcal G}$ is known as
\begin{equation}\label{eq:def:leqepsilon}
\norm{A_{\Delta\mathcal G}} \leq \epsilon
\end{equation}
for a positive constant~$\epsilon$. 

We consider the following standard epidemiological problem (see \cite{Preciado2014} for the case where no uncertainty exists in the underlying network). We assume that we can distribute within the network vaccines for reducing the infection rates of individuals, and antidotes for increasing their recovery rates. Let us suppose that the infection and recovery rates can be tuned within the intervals
\begin{equation}\label{eq:bosConsts}
\ubar{\beta} \leq \beta_i \leq \bar{\beta}, \ \ubar{\delta} \leq \delta_i \leq \bar{\delta}. 
\end{equation}
Let $f(\beta_i)$ denote the cost for setting the infection rate of node~$i$ to $\beta_i$. Likewise, let $g(\delta_i)$ denote the cost for setting the recovery rate of node~$i$. Then, the total cost $L$ for achieving a set of infection and recovery rates~$(\beta_1, \dotsc, \beta_N,  \delta_1, \dotsc, \delta_N)$ equals 
\begin{equation}\label{eq:totalCost}
L = \sum_{i=1}^N(f(\beta_i) + g(\delta_i)). 
\end{equation}
Through the resource distribution, we aim for increasing the exponential decay rate of the epidemic process defined by 
\begin{equation*}
\rho = -\sup_{\mathcal V_0\subset \mathcal V}\limsup_{t\to\infty}\frac{\log \sum_{i=1}^N p_i(t)}{t}, 
\end{equation*}
where $\mathcal V_0$ denotes the set of initially infected nodes and~$p_i(t)$ denotes the probability that node~$i$ is infected at time~$t$. We can now state the resource distribution problem studied in this section. 

\begin{problem}\label{prb:}
Let a minimum required exponential decay rate, denoted by~$\gamma > 0$, be given. Find the set of infection rates~$\{\beta_i\}_{i=1}^N$ and recovery rates~$\{\delta_i\}_{i=1}^N$ that minimizes the total cost~$L$ given by \eqref{eq:totalCost}, while satisfying the following robust stability condition
\begin{equation}\label{eq:wanttoachieve}
\inf_{\Delta\mathcal G\colon \norm{A_{\Delta \mathcal G}}\leq \epsilon} \rho > \gamma. 
\end{equation}
\afterequation
\end{problem}

The computation of the exponential decay rate~$\rho$ is very hard for contact networks of large size because of the huge size of the state space of the SIS model (as a Markov process). A popular approach to simplify the analysis of this type of Markov processes is to consider upper-bounding linear models (see, e.g.,~\cite{Preciado2014}), from which we obtain 
\begin{equation*}
\rho \geq -\lambda_{\max}(\diag(\beta) A_{\mathcal G}- \diag(\delta)). 
\end{equation*}
Therefore, to satisfy the robust stability condition \eqref{eq:wanttoachieve}, it is sufficient to achieve that
\begin{equation}\label{eq:SIScondition}
\sup_{\Delta\mathcal G\colon \norm{A_{\Delta \mathcal G}}\leq \epsilon} \lambda_{\max}\left(\diag(\beta) (A_{\bar{\mathcal G}} + A_{\Delta \mathcal G})- \diag(\delta)\right) < -\gamma. 
\end{equation}
We use this fact to reduce Problem~\ref{prb:} into a robust stabilization problem of the form~\eqref{eq:strucutredStabilization}. Let us introduce the vectorial parameter
\begin{equation}\label{eq:def:theta}
\theta = [\beta_1, \, \dotsc, \beta_N,\, \delta_1, \, \dotsc, \,  \delta_N]^\top. 
\end{equation}
Define 
\begin{equation}\label{eq:def:Atheta}
A(\theta) = \diag(\beta) A_{\mathcal G}- \diag(\delta), 
\end{equation}
$B(\theta) = \diag(\beta)$, and~$C(\theta) =I$. Then, we can rewrite the requirement~\eqref{eq:SIScondition} as \eqref{eq:strucutredStabilization:eta}. Therefore, Problem~\ref{prb:} reduces to the robust stabilization problem~\eqref{eq:strucutredStabilization} studied in Section~\ref{sec:robust}.

\begin{figure}
\centering
\includegraphics[width=.95\linewidth,trim={0cm .1cm 0cm 0cm},clip]{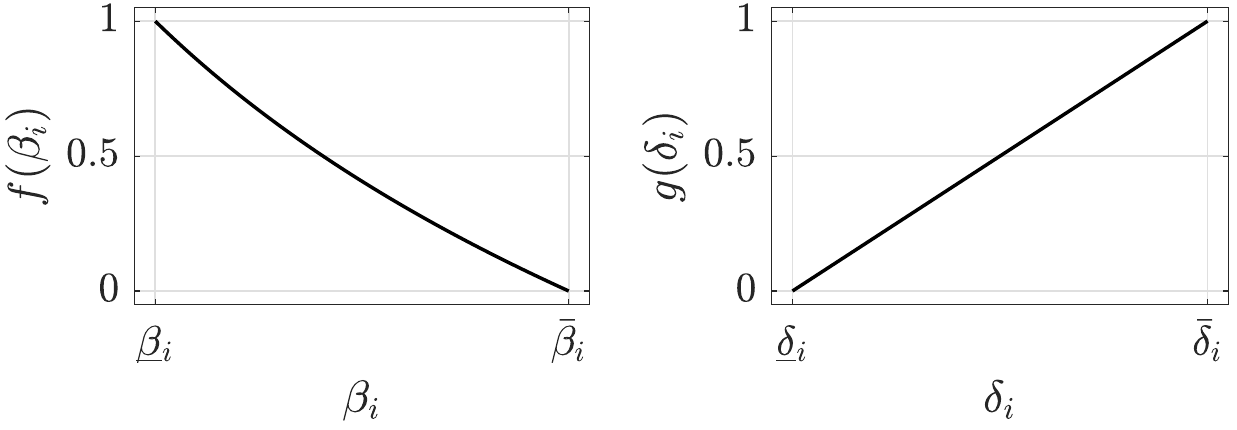}
\caption{Cost functions for infection and recovery rates.}
\label{fig:costFunctions}
%
%
\vspace{6mm}
\centering
\includegraphics[height=.6\linewidth,trim={0cm 0cm 0cm .5cm},clip]{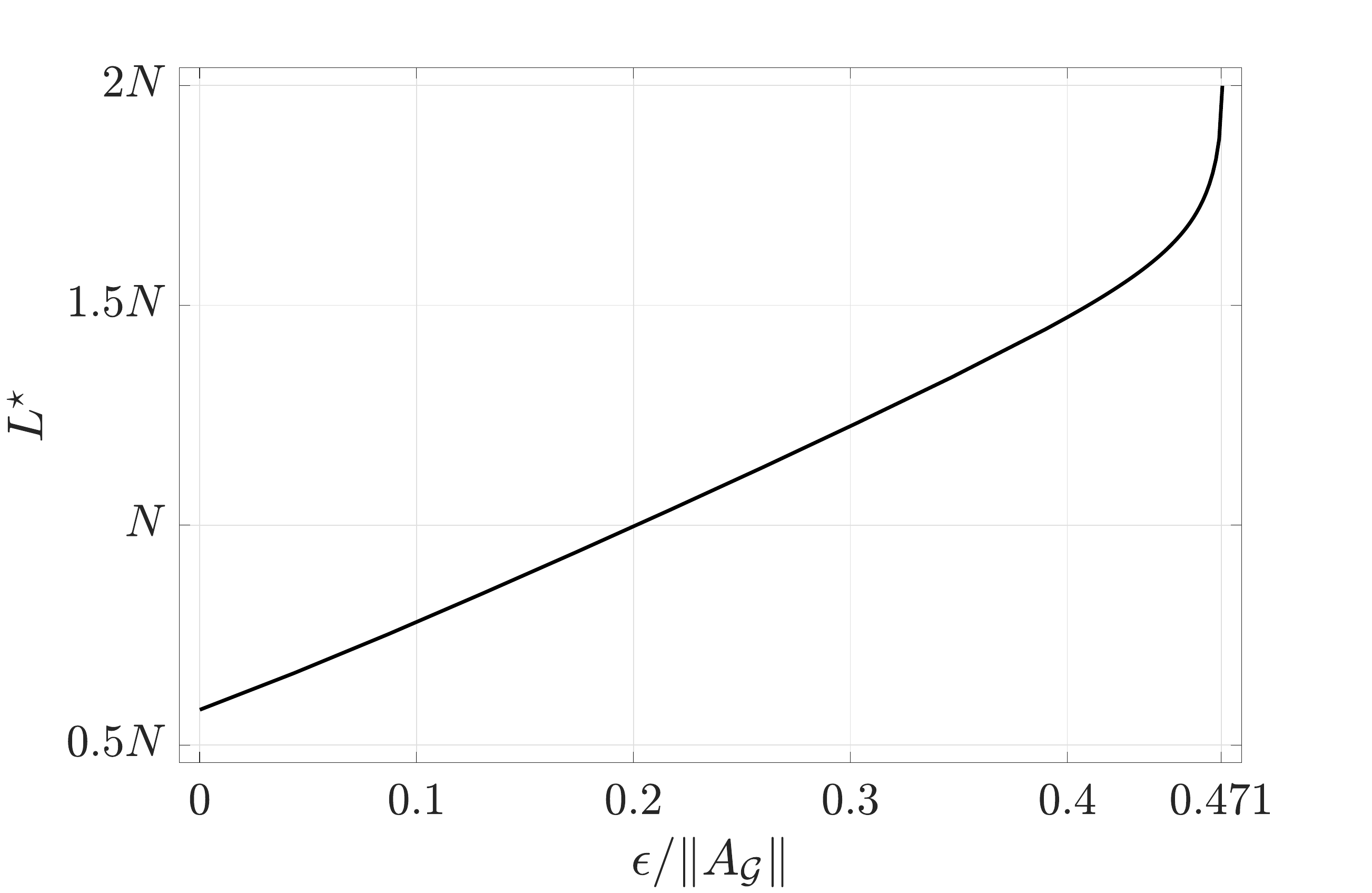}
\caption{The amount of optimal total investments~$L^\star$ versus the relative maximum size~$\epsilon/\norm{A_{\mathcal G}}$ of the additive uncertainty~$\norm{A_{\Delta \mathcal G}}$}
\label{fig:errorBudget1}
\end{figure}

In order to apply Theorem~\ref{thm:mu} for solving Problem~\ref{prb:} via geometric programming, we need to confirm that Assumptions~\ref{asm:posy} and~\ref{asm:costPosy} hold true. It is easy to see that Assumption~\ref{asm:posy} is satisfied because, for the diagonal matrix~
\begin{equation}\label{eq:R=diagdelta}
R(\theta) = \diag(\delta)
\end{equation}
with monomial diagonals, each entry of the matrix function~$\tilde A(\theta) = A(\theta) + R(\theta) = \diag(\beta) A_{\mathcal G}$ is either a posynomial with respect to the variables in~\eqref{eq:def:theta} or zero. To guarantee that Assumption~\ref{asm:costPosy} holds true, let us use the following cost functions similar to the ones used in~\cite{Preciado2014}:
\begin{equation}\label{eq:costdeff}
f(\beta_i) = \frac{\beta_i^{-p} - \bar \beta^{-p}}{\ubar \beta^{-p} - \bar \beta^{-p}},\ g(\delta_i) = \frac{\delta_i^q - \ubar \delta^q}{\bar \delta^q - \ubar \delta^q}, 
\end{equation}
where $p>0$ and~$q>0$ are constants to tune the shape of the cost functions. Notice that the cost function~$f$ is normalized as $f(\ubar \beta) = 1$ and~$f(\bar \beta) = 0$. This indicates that $\bar \beta$ is the nominal infection rate of nodes, and that a unit investment improves the nominal rate to the minimum possible infection rate~$\ubar \beta$. The same interpretation applies to the cost function~$g$ for recovery rates. When the above cost functions are used, the total cost~$L$ in~\eqref{eq:totalCost} satisfies Assumption~\ref{asm:costPosy}.\ref{ams:item:L}) with the constant
\begin{equation*}
L_0 = N \left(
\frac{\bar \beta^{-p}}{\ubar \beta^{-p} - \bar \beta^{-p}}+\frac{\ubar \delta^q}{\bar \delta^q - \ubar \delta^q} 
\right). 
\end{equation*}
Also, the box constraints~\eqref{eq:bosConsts} can be easily converted to constraints in terms of posynomials. Therefore, Assumption~\ref{asm:costPosy} is satisfied as well.

In this numerical simulation, we let the nominal network~$\bar{\mathcal G}$ be a human social network having $N=247$ nodes with its adjacency matrix having spectral radius~$13.53$. Suppose that $\ubar \beta = 0.1$, $\bar \beta = 0.2$, $\ubar \delta = 1$, and~$\bar \delta = 2$. The exponents $p, q$ in the cost functions~\eqref{eq:costdeff} are chosen as $p=0.1$ and~$q=1$. The graphs of the corresponding cost functions are shown in Fig.~\ref{fig:costFunctions}. We require that the exponential decay rate of the SIS model is at least $\gamma = 0.01$ for any additive uncertainty $A_{\Delta \mathcal G}$ satisfying inequality~\eqref{eq:def:leqepsilon}.

\begin{figure}
\centering
\subfloat[$\epsilon = 0$]{\includegraphics[height=.6\linewidth,trim={0cm 0cm 2cm .75cm},clip]{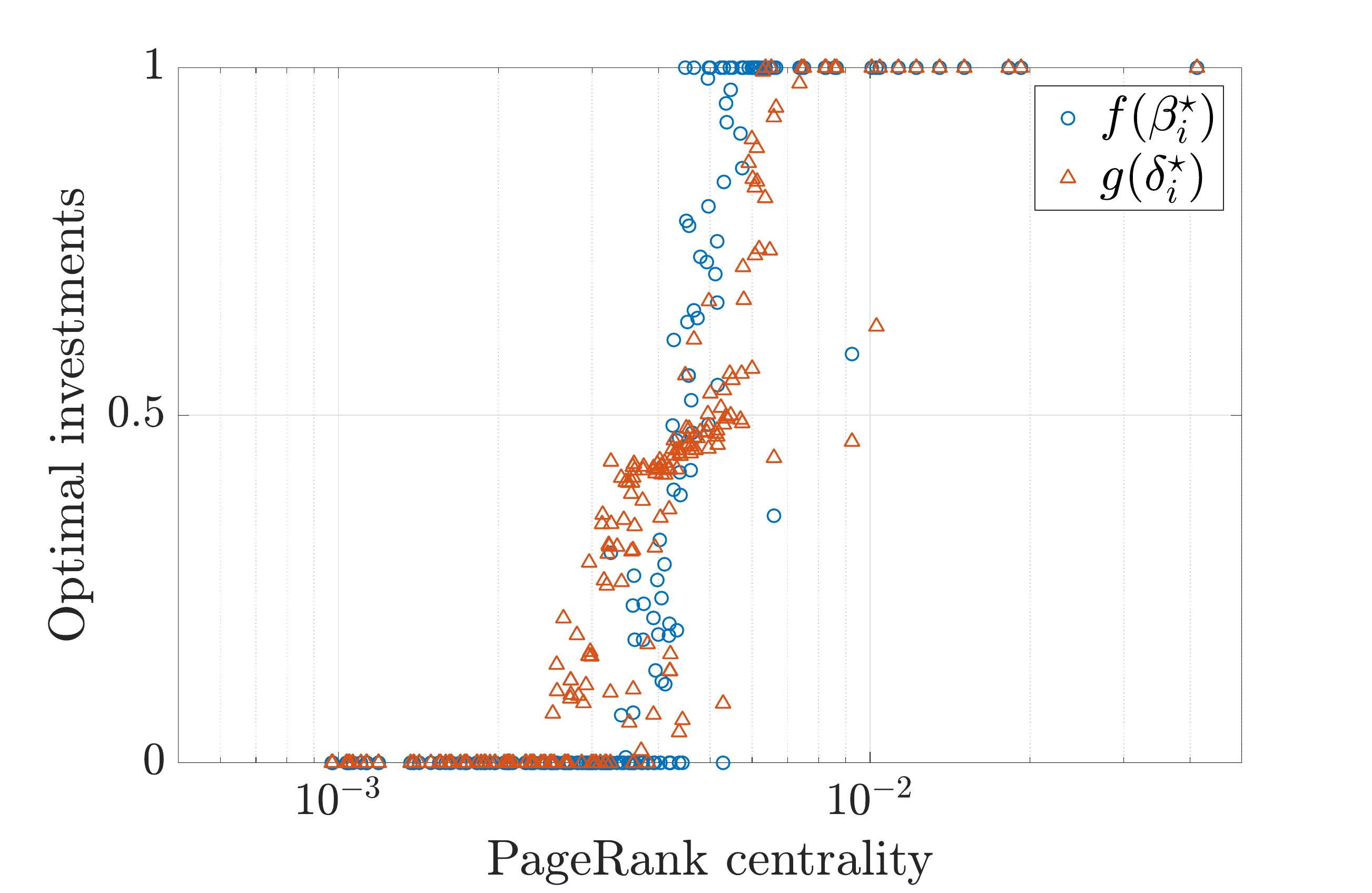}
\label{fig:investmentFirst}}
\\
\subfloat[{$\epsilon = 0.4\cdot \norm{A_{\mathcal G}}$}]{\includegraphics[height=.6\linewidth,trim={0cm 0cm 2cm .75cm},clip]{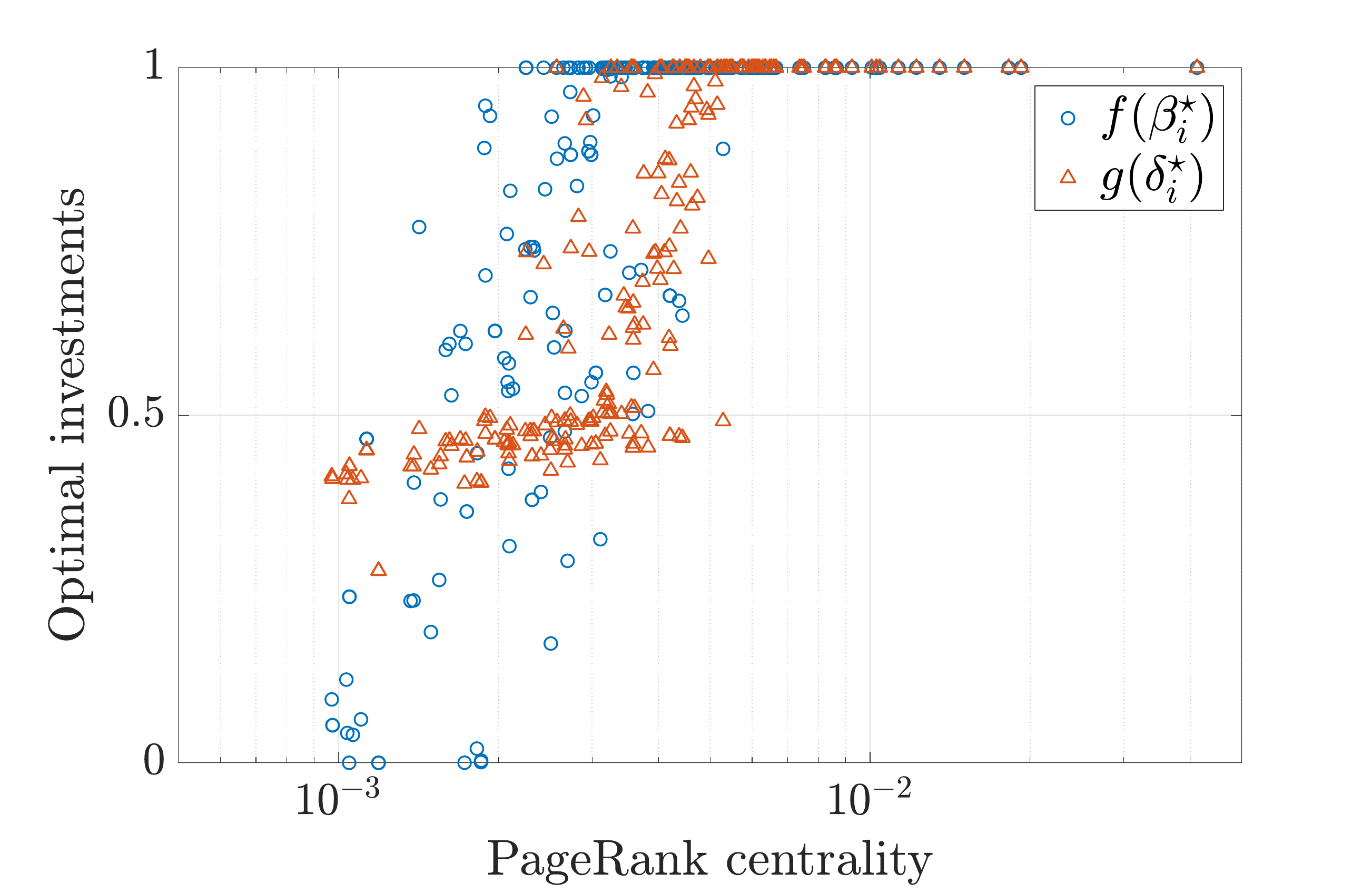}\label{fig:investmentLast}
}
\caption{Optimal investments $f(\beta_i^\star), g(\delta_i^\star)$ on individual nodes versus the PageRank centralities of the nodes.}
\label{fig:investment}
\end{figure}

We first use Corollary~\ref{cor:} and identify the maximum allowable size of the uncertainty as $\epsilon^\star = 0.471\cdot \norm{A_{\mathcal G}}$. Then, for various values of~$\epsilon$ in the interval $[0, \epsilon^\star]$, we solve the resource distribution problem to find the optimized infection rates~$\beta_i^\star$ and recovery rates~$\delta_i^\star$ by geometric programming. In Fig.~\ref{fig:errorBudget1}, we show how the optimal total cost, denoted by~$L^\star$, depends on the size of the uncertainty~$\epsilon$. We then investigate how a particular value of the uncertainty size~$\epsilon$ affects the way in which medical resources are distributed over the complex network. In Fig.~\ref{fig:investment}, we show the amount of resources spent on improving the infection and recovery rates of individual nodes versus the PageRank~\cite{Page1998} of the nodes in the nominal network, for the cases of~$\epsilon=0$ (Fig.~\ref{fig:investmentFirst}) and~$\epsilon = 0.4\cdot \norm{A_{\mathcal G}}$ (Fig.~\ref{fig:investmentLast}), respectively. When no uncertainty is expected ($\epsilon=0$), we find several nodes not receiving investments on their recovery rates. This trend drastically disappears as we increase the size of the uncertainty to $\epsilon = 0.4\cdot \norm{A_{\mathcal G}}$, in which case all nodes receive at least one-fourth unit of investments on their recovery rates. This observation indicates the importance of correctly identifying the connectivity structure of complex networks for effective distribution of medical resources.

{\begin{remark}\label{rmk:diagonal}
The diagonal matrix~$R(\theta)$ in \eqref{eq:R=diagdelta} cannot be written in the form~\eqref{eq:Rfurtherassumption}. This fact seems to prevent us from using Theorems~\ref{thm:h2cont}, \ref{thm:Hankel}, and \ref{thm:Schatten} for optimizing the $H^2$, Hankel, and Schatten-$p$ norms of the epidemic dynamics described by the SIS model. We can, however, avoid this problem by using a different parametrization of system matrices. Let us parametrize the recovery rate~$\delta_i$ as 
\begin{equation*}
\delta_i = \bar \delta + 1 - \frac{1}{\delta_i^{(c)}},\quad  
\frac{1}{\bar \delta - \ubar{\delta}+1}\leq \delta_i^{(c)}\leq 1, 
\end{equation*}
where $\delta_i^{(c)}$ is an auxiliary positive variable. 
Let us introduce the notation $\diag(1/\delta^{(c)}) = \diag(1/\delta^{(c)}_1, \cdots, 1/\delta^{(c)}_N)$. Then, we can rewrite the matrix~\eqref{eq:def:Atheta} as $A(\theta) = \tilde A(\theta) - R(\theta)$ for $\tilde A(\theta) = \diag(\beta)A_{\mathcal G} + \diag(1/\delta^{(c)})$ and $R(\theta) = - (\bar \delta + 1)I$. This $R(\theta)$ now is of the form~\eqref{eq:Rfurtherassumption}. Furthermore, each entry of the matrix~$\tilde A(\theta)$ is indeed a posynomial of the variables, as desired. 
\end{remark}}

\section{Conclusions and discussion} \label{sec:conclusion}

In this paper, we have presented geometric programming-based frameworks for the parameter tuning problem of positive linear systems constrained by a parameter tuning cost as well as system norms or stability properties. We have considered the following standard system norms; the $H^2$ norm, $H^\infty$ norm, Hankel norm, and Schatten $p$-norm. We have also shown that the robust stabilization problem under structured uncertainties, as well as a mixed-constraint parameter tuning problem for delayed positive linear systems can be numerically efficiently solved. We have illustrated the effectiveness of our theoretical results via numerical simulations on dynamical buffer networks and networked epidemic spreading processes. 

There are several research directions that should be further pursued. One such direction is the synthesis of switched positive linear systems~\cite{Gurvits2007,Blanchini2012,Zhao2012}. In particular, it is of theoretical interest to investigate if we can utilize linear programming-based results for the analysis of positive Markov jump linear systems (see, e.g.,~\cite{Bolzern2014}) to obtain geometric programs for synthesis problems. Another research direction of interest is the synthesis of cone-preserving linear systems. It has been found in the literature~\cite{Tanaka2013,Shen2017a,Shen2017b} that linear systems leaving a cone invariant share several interesting properties with positive linear systems. In this direction, it is left as an open problem to examine if the current geometric programming-based approach can be applied to cone-preserving linear systems.





\begin{IEEEbiography}[{\includegraphics[width=1in,height=1.25in,clip,keepaspectratio]{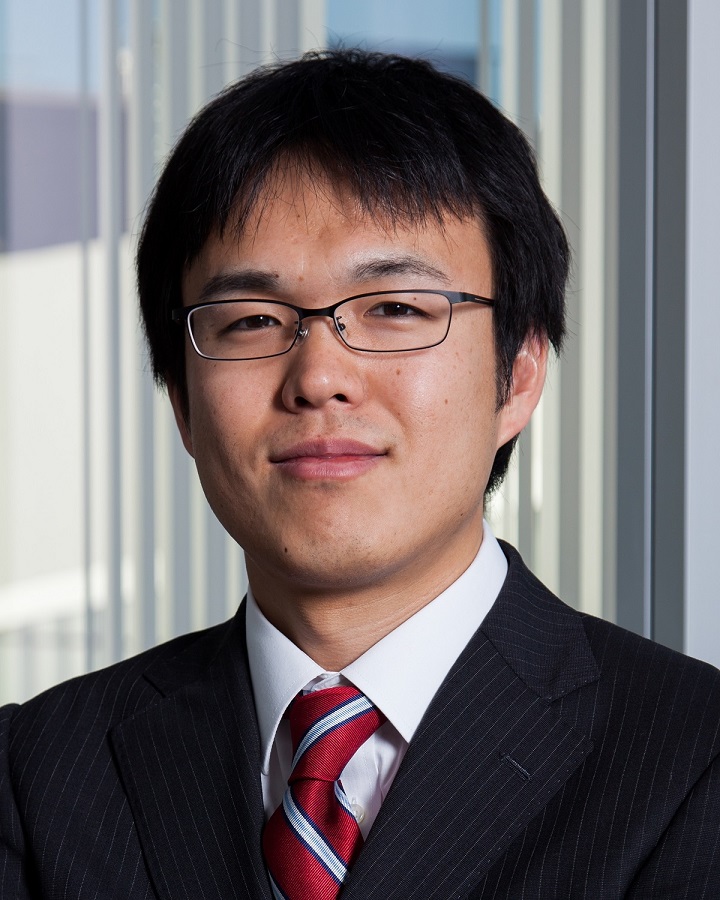}}]{Masaki Ogura}
Masaki Ogura is an Associate Professor in the Graduate School of Information Science and Technology at Osaka University, Japan. He received his M.Sc. degree in Informatics from Kyoto University in 2009, and his Ph.D.~in Mathematics from Texas Tech University in 2014. From 2014 to 2017, he was a Postdoctoral Researcher at the University of Pennsylvania. From 2017 to 2019, he was an Assistant Professor at the Nara Institute of Science and Technology, Japan. His research interests include network science, dynamical systems, and stochastic processes with applications in networked epidemiology, design engineering, and biological physics. He was a runner-up of the 2019 Best Paper Award by the IEEE Transactions on Network Science and Engineering and a recipient of the 2012 SICE Best Paper Award. 
\end{IEEEbiography}

\begin{IEEEbiography}[{\includegraphics[width=1in,height=1.25in,clip,keepaspectratio]{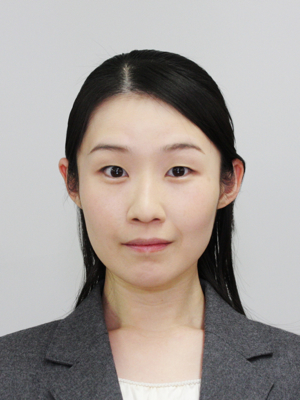}}]{Masako Kishida}
Dr. Kishida received her Ph.D. from the University of Illinois at Urbana-Champaign in 2010. After holding appointments at universities in the U.S.A., Japan, New Zealand, and Germany, she became Associate Professor at the National Institute of Informatics, Tokyo Japan, in 2016. She received Humboldt Research Fellowship from the Alexander von Humboldt Foundation in 2015 and Telecom System Technology Award from the Telecommunications Advancement Foundation in 2019. She is a senior member of IEEE since 2018.
\end{IEEEbiography}


\begin{IEEEbiography}[{\includegraphics[width=1in,height=1.25in,clip,keepaspectratio]{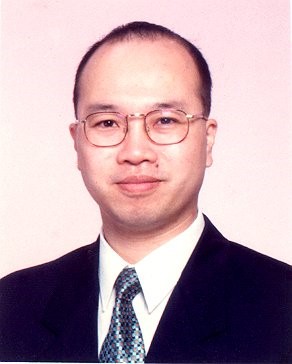}}]{James Lam}
Professor J.~Lam received a BSc (1st Hons.) degree in Mechanical Engineering from the University of Manchester, and was awarded the Ashbury Scholarship, the A.H. Gibson Prize, and the H. Wright Baker Prize for his academic performance. He obtained the MPhil and PhD degrees from the University of Cambridge. He is a Croucher Scholar, Croucher Fellow, and Distinguished Visiting Fellow of the Royal Academy of Engineering. Prior to joining the University of Hong Kong in 1993 where he is now Chair Professor of Control Engineering, he was a lecturer at the City University of Hong Kong and the University of Melbourne. 

Professor Lam is a Chartered Mathematician, Chartered Scientist, Chartered Engineer, Fellow of Institute of Electrical and Electronic Engineers, Fellow of Institution of Engineering and Technology, Fellow of Institute of Mathematics and Its Applications, Fellow of Institution of Mechanical Engineers, and Fellow of Hong Kong Institution of Engineers. He is Editor-in-Chief of IET Control Theory and Applications and Journal of The Franklin Institute, Subject Editor of Journal of Sound and Vibration, Editor of Asian Journal of Control, Senior Editor of Cogent Engineering, Associate Editor of Automatica, International Journal of Systems Science, Multidimensional Systems and Signal Processing, and Proc. IMechE Part I: Journal of Systems and Control Engineering. He is a member of the Engineering Panel (Joint Research Scheme), Research Grant Council, HKSAR. His research interests include model reduction, robust synthesis, delay, singular systems, stochastic systems, multidimensional systems, positive systems, networked control systems and vibration control. He is a Highly Cited Researcher in Engineering (2014, 2015, 2016, 2017, 2018{, 2019}) and Computer Science (2015).
\end{IEEEbiography}

\end{document}